\documentclass{conm-p-l}
\usepackage{epsf}
\usepackage[all]{xy}

\newtheorem{theorem}{Theorem}[section]
\newtheorem{lemma}[theorem]{Lemma}

\theoremstyle{definition}

\newtheorem{proposition}[theorem]{Proposition}

\theoremstyle{remark}

\numberwithin{equation}{section}

\newcommand{\abs}[1]{\lvert#1\rvert}

\newcommand{\internalcomment}[1]{}

\newcommand{\Z}{\mathbf{Z}}

\newcommand{\bP}{\mathbf{P}}

\newcommand{\ko}{\: , \;}
\newcommand{\ie}{{\em i.~e.}\ }
\newcommand{\cf}{{\em cf.}\ }
\newcommand{\eg}{{\em e.~g.}\ }

\newcommand{\ol}{\overline}

\renewcommand{\tilde}[1]{\widetilde{#1}}

\newcommand{\ra}{\rightarrow}
\newcommand{\la}{\leftarrow}

\newcommand{\arr}[1]{\stackrel{#1}{\rightarrow}}

\newcommand{\iso}{\stackrel{_\sim}{\rightarrow}}
\newcommand{\liso}{\stackrel{_\sim}{\leftarrow}}

\newcommand{\opname}[1]{\operatorname{\mathsf{#1}}}

\newcommand{\sets}{\opname{Sets}}

\newcommand{\aqu}{{\opname{aQu}\nolimits}}

\newcommand{\aco}{{\opname{aCo}\nolimits}}
\newcommand{\adgco}{{\opname{adgCo}\nolimits}}
\newcommand{\agco}{{\opname{agCo}\nolimits}}
\newcommand{\For}{\opname{For}}

\newcommand{\alg}{\opname{Alg}}
\newcommand{\cog}{\opname{Cog}}

\newcommand{\Ho}{\opname{Ho}}
\renewcommand{\mod}{\opname{mod}\nolimits}
\newcommand{\Mod}{\opname{Mod}\nolimits}
\newcommand{\Comc}{\opname{Comc}\nolimits}
\newcommand{\Inj}{\opname{Inj}\nolimits}

\newcommand{\Grmod}{\opname{Grmod}\nolimits}

\newcommand{\coh}{\opname{coh}\nolimits}
\newcommand{\obj}{\opname{obj}\nolimits}

\newcommand{\id}{\mathbf{1}}

\newcommand{\ten}{\otimes}

\newcommand{\tp}[1]{^{\ten #1}}

\newcommand{\im}{\opname{im}\nolimits}
\renewcommand{\ker}{\opname{ker}\nolimits}

\newcommand{\HH}[1]{{HH}^{i}\,}

\newcommand{\ca}{{\mathcal A}}
\newcommand{\cb}{{\mathcal B}}
\newcommand{\cc}{{\mathcal C}}

\newcommand{\cF}{{\mathcal F}}

\newcommand{\ct}{{\mathcal T}}

\newcommand{\cw}{{\mathcal W}}

\newcommand{\eps}{\varepsilon}

\renewcommand{\phi}{\varphi}

\newcommand{\End}{\opname{End}}
\newcommand{\Hom}{\opname{Hom}}
\newcommand{\Homdot}{\opname{Hom^\bullet}}
\newcommand{\Homq}{\underline{Hom}_q}
\newcommand{\Homc}{\underline{Hom}_c}
\newcommand{\Funinf}{\opname{Fun}_\infty}
\newcommand{\Ext}{\opname{Ext}}

\newcommand{\tensinf}{\overset{\infty}{\ten}}
\newcommand{\Tor}{\opname{Tor}}
\newcommand{\Cotor}{\opname{Cotor}}

\newcommand{\Tw}{\opname{Tw}}

\newcommand{\fun}{\opname{fun}}

\newcommand{\per}{\opname{per}}
\newcommand{\der}{\opname{D}}

\newcommand{\fib}{\cF ib}
\newcommand{\cof}{\cc of}
\newcommand{\centeps}[1]{\begin{array}{c} \epsfbox{#1} \end{array}}

\begin{document}

\title{A-infinity algebras, modules and functor categories}

\author{Bernhard Keller}
\address{UFR de Math\'ematiques, Universit\'e Paris 7, Institut de Math\'ematiques, UMR 7586 du CNRS,
Case 7012, 2 place Jussieu, 75251 Paris cedex 05, France}

\email{keller@math.jussieu.fr}


\subjclass{18E30, 16D90, 18G40, 18G10, 55U35}
\date{September 7, 2005}


\keywords{A-infinity structure, strong homotopy algebra, Quillen
model category, triangulated category}

\begin{abstract}
In this survey, we first present basic facts on A-infinity
algebras and modules including their use in describing triangulated
categories. Then we describe the Quillen model approach to A-infinity
structures following K.~Lef\`evre's thesis.
Finally, starting from an idea of V.~Lyubashenko's,
which is being fully developped in \cite{BespalovLyubashenkoManzyuk06}, we
give a conceptual construction of A-infinity functor
categories using a suitable closed monoidal category of
cocategories. In particular, this yields a natural construction
of the bialgebra structure on the bar construction of the
Hochschild complex of an associative algebra.
\end{abstract}

\maketitle

\section{Introduction}

A-infinity spaces and
A-infinity algebras were invented at the beginning of the
sixties by Stasheff \cite{Stasheff63}. In the seventies and the eighties,
they were developped further by Smirnov \cite{Smirnov80},
Kadeishvili \cite{Kadeishvili80}, Prout\'e \cite{Proute84},
Huebschmann \cite{Huebschmann86}, \cite{Huebschmann89},
\ldots
especially with a view towards applications in topology.
At the beginning of the nineties, the relevance of
A-infinity structures in geometry and physics became
apparent through the work of Getzler-Jones \cite{GetzlerJones90},
Stasheff \cite{Stasheff92}, Fukaya \cite{Fukaya93}, Kontsevich
\cite{Kontsevich94}, \ldots,
and later Kontsevich-Soibelman \cite{KontsevichSoibelman01},
Seidel \cite{Seidel02}, \ldots. We refer to \cite{Keller01},
\cite{Keller02} for a more detailed introduction with numerous references.
The present survey consists of three parts:
\begin{itemize}
\item[1.]  In sections 1--3, we will define $A_\infty$-algebras
and examine their basic properties. Then we will define
$A_\infty$-modules, the derived category and conclude with
the description of triangulated categories via $A_\infty$-algebras.
\item[2.] In section 4, we will present the interpretation of
$A_\infty$-algebras (resp. $A_\infty$-modules) as the fibrant objects
in the model category of certain differential graded coalgebras
(resp. comodules) following Lef\`evre \cite{Lefevre03}.
\item[3.] In section 5, starting from an idea of
V.~Lyubashenko's \cite{Lyubashenko03}
we will give a conceptual construction
of $A_\infty$-functor categories and use it to construct the
canonical bialgebra structure on the cobar construction of the Hochschild
complex of an associative algebra.
\end{itemize}

\noindent {\bf Acknowledgments.} The first two parts of this
survey are based on lectures given by the author at the Workshop
on Derived Categories, Quivers and Strings (Edinburgh, August
2004) and at the Conference on Topology (Canberra, July 2003). The
author is grateful to the organizers of both events and in
particular to Alastair King and Amnon Neeman for their interest
and encouragment. He thanks Vladimir Lyubashenko for stimulating
discussions on the third part and Pedro Nicol\'as and Oleksandr
Manzyuk for carefully reading the manuscript, pointing out
and correcting errors in previous versions of
section~\ref{ss:Augmented-k-quivers}. He is grateful to the
organizers of the ICRA XIII for including this material in the
proceedings of the workshop.

\section{A-infinity algebras}

\subsection{Notations} We will follow Fukaya's sign and
degree conventions. For this, we need to introduce some notation:

Let $k$ be a field. If $V$ is a graded vector space, \ie
$$
V = \bigoplus_{p\in\Z} V^p \ko
$$
we denote by $SV$ or $V[1]$ the graded space with $(SV)^p=V^{p+1}$
for all $p\in\Z$. We call $SV$ the {\em suspension} or the
{\em shift} of $V$.

If $f: V \to V'$ and $g: W \to W'$ are homogeneous maps between
graded spaces, their {\em tensor product}
$$
f \ten g : V \ten W \to V'\ten W'
$$
is defined by
$$
(f \ten g)(v \ten w)= (-1)^{\abs{g}\, \abs{v}} f(v) \ten g(w)
$$
for all homogeneous elements $v\in V$ and $w\in W$.

If $V$ and $V'$ are complexes, \ie endowed with differentials
$d$ homogeneous of degree $1$ and square $0$, we put $d_{SV}=-d_V$
and, for a homogeneous map $f:V \to V'$,
$$
d(f)= d_{V'} \circ f - (-1)^{\abs{f}} f \circ d_V.
$$
Thus, $f$ is a morphism of complexes iff $d(f)=0$ and
two morphisms of complexes $f$ and $f'$ are homotopic iff there is
a morphism of graded spaces $h$ such that $f'=f+d(h)$.
We will use that homotopic morphisms induce the same map in homology.

\subsection{A-infinity algebras}
An {\em $A_\infty$-algebra} is a graded space $A$ endowed with
maps
$$
b_n : (SA)^{\ten n} \to SA
$$
defined for $n\geq 1$, homogeneous of degree $1$ and such that,
for all $n\geq 1$, we have
\begin{equation} \label{ainfeq}
\sum_{i+j+l=n} b_{i+1+l} \circ ( \id\tp{i}\ten b_j \ten \id\tp{l}) = 0
\end{equation}
as maps from $(SA)\tp{p}$ to $SA$. Here, the symbol $\id$ denotes the
identity map of $SA$. We visualize $b_n$ either as a planar
tree with $n$ leaves and one root or as a halfdisk whose
upper arc is divided into segments, each of which symbolizes
an \lq input\rq, and whose base segment symbolizes the \lq output\rq\
of the operation.
\medskip
$$
\centeps{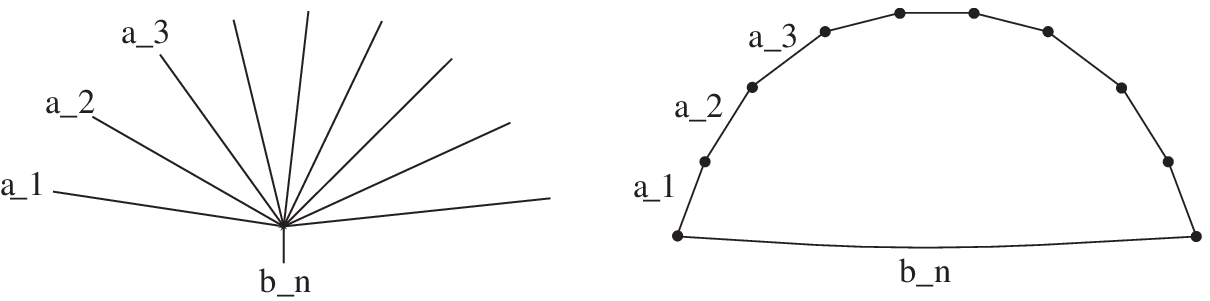}
$$
Using this last representation, the defining
identity~(\ref{ainfeq}) is depicted as follows:
\medskip
\[
\sum \pm \centeps{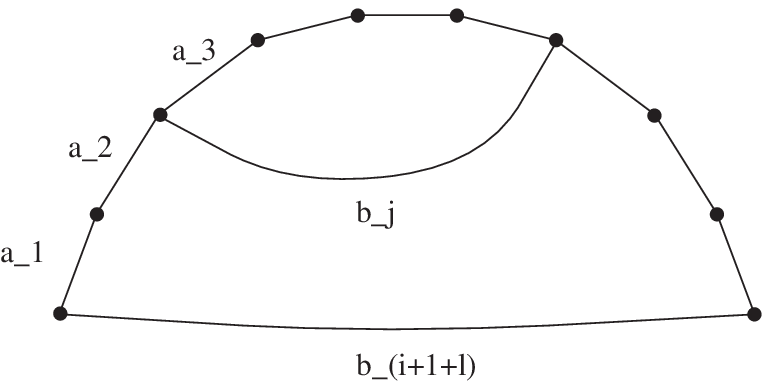}  = 0.
\]
Let us analyze the defining identity for small values of $n$:
For $n=1$, it states that $b_1^2=0$, so that $(SA, b_1)$
is a complex. We also make $A$ into a complex by endowing
it with the shifted differential
$$
m_1=-b_1.
$$
For $n=2$, the defining identity becomes
$$
b_1 b_2 + b_2 (b_1 \ten \id + \id \ten b_1) =0.
$$
Note that $b_1 \ten \id+ \id \ten b_1$ is the differential
of $SA\ten SA$ so that we obtain $d(b_2)=0$, which means that
$b_2 : SA\ten SA \to SA$ is a morphism of complexes.

For $n=3$, the identity~(\ref{ainfeq}) becomes
$$
b_2 \, (b_2 \ten \id + \id \ten b_2) + b_1 b_3 + b_3 (b_1 \ten \id \ten \id
+ \id\ten b_1 \ten \id + \id \ten \id \ten b_1) =0.
$$
Here the second summand is $d(b_3)$ whereas the first is,
up to a sign, the associator for the binary operation $b_2$.
If we define $m_2: A \ten A \to A$ by
$$
m_2(x,y)= (-1)^{\abs{x}} b_2(x,y),
$$
then we obtain that $m_2$ is associative up to a homotopy
given by $b_3$.

For each $n>3$, the identity~(\ref{ainfeq}) states that the
maps $b_2, \ldots, b_{n-1}$ satisfy a certain quadratic
identity up to a homotopy given by $b_n$. In this sense,
an $A_\infty$-algebra is an algebra associative up to a
given system of higher homotopies.

It is a direct consequence of the definition, that if
$b_n$ vanishes for each $n\geq 3$, then $(A, m_1, m_2)$
is a differential graded (=dg) algebra, \ie $m_2$ is
associative and $m_1$ a differential compatible with
$m_2$ through the (graded) Leibniz rule. Conversely,
each dg algebra gives rise to an $A_\infty$-algebra
with vanishing $b_n$, $n\geq 3$.

In particular, each ordinary associative algebra can
be viewed as an $A_\infty$-algebra concentrated in degree
$0$, and conversely.

\subsection{Examples via deformations} Following an idea
of Penkava-Schwarz \cite{PenkavaSchwarz95},
let us exhibit a large class of easily
constructed but non trivial examples of $A_\infty$-algebras.
Let $B$ be an ordinary algebra and $N\geq 1$ an integer.
Let $\eps$ be an indeterminate of degree $2-N$. We first
endow the graded space $A=B[\eps]/(\eps^2)$ with the
trivial $A_\infty$-structure given by the map $b_2$
induced by the multiplication of $B$ and the maps
$b_n=0$ for all $n\neq 2$. Now let
$$
c: B\tp{N} \to B
$$
be any linear map. Define deformed multiplications
$$
b'_n = \left\{ \begin{array}{ll} b_n & n \neq N \\ b_N + \eps c & n=N .
\end{array} \right.
$$
Then it is easy to see that {\em $A$ endowed with the $b'_n$ is
an $A_\infty$-algebra iff $c$ is a Hochschild cocycle for $B$}.

\subsection{Weak A-infinity algebras} A {\em weak
$A_\infty$-algebra} is a graded space $A$ endowed with
maps $b_0 : k \to SA$ and $b_n$, $n\geq 0$, such that
the identity~\ref{ainfeq} holds for all $n\geq 0$.
The preceding example then naturally extends to the case
where $N=0$, where we start from a Hochschild
$0$-cocycle, \ie a central element $c$ of $B$.
In general, in a weak $A_\infty$-algebra, we have
$$
b_1^2 = - b_2(b_0 \ten \id + \id \ten b_0) \neq 0
$$
so that the homology with respect to $b_1$ is
no longer defined and the above remarks no longer apply.
Little is known about weak $A_\infty$-algebras in general,
but they do appear in nature as deformations.

\subsection{Morphisms and quasi-isomorphisms}
\label{ss:Morphisms-and-quasi-isomorphisms}
A {\em morphism of $A_\infty$-algebras $f:A \to B$} is given
by maps
$$
f_n : (SA)\tp{n} \to SB \ko n\geq 1 \ko
$$
homogeneous of degree $0$ such that, for all $n\geq 1$,
we have
$$
\sum_{i+j+l=n} f_{i+1+l} \circ ( \id\tp{i} \ten b_j \ten \id\tp{l})
= \sum_{i_1+ \cdots i_s=n} b_s \circ (f_{i_1} \ten \cdots \ten f_{i_s}).
$$
By looking at this equation for $n=1$ and $n=2$ we see that
$f_1$ then induces a morphism of complexes from $(A,m_1)$ to
$(B,m_1)$ which is compatible with $m_2$ up to an homotopy
given by $f_2$. In particular, $f_1$ induces an {\em algebra
morphism}
$$
H^*A \to H^* B.
$$
By definition, $f$ is an $A_\infty$-quasi-isomorphism if
$f_1$ is a quasi-isomorphism (\ie induces an isomorphism in
homology).

The {\em composition $f\circ g$} of two morphisms is given by
$$
(f\circ g)_n = \sum_{i_1+\cdots + i_s =n} f_{s} \circ (g_{i_1}\ten
\cdots g_{i_s}).
$$
The {\em identical morphism of $SA$} is given by $f_1=\id$ and
$f_n=0$ for all $n\geq 2$.

It is easy to see that we do obtain a category. It contains the
category of dg algebras and their morphisms as a non-full subcategory.

\begin{proposition} (\cf section~\ref{link-Ainf-modules})
For each $A_\infty$-algebra $A$, there is a universal $A_\infty$-algebra
morphism $\phi : A \to U(A)$ to a dg algebra $U(A)$. Moreover, $\phi$ is
an $A_\infty$-quasi-isomorphism.
\end{proposition}

The universal property means that for each dg algebra $B$, each
$A_\infty$-morphism $f: A \to B$ factors as $f=g\circ \phi$ for
a unique morphism of dg algebras $g: U(A) \to B$. The proposition
tells us that, up to $A_\infty$-quasi-isomorphism, $A_\infty$-algebras
are quite similar to dg algebras. However, in other respects, they
are radically different from dg algebras, as the following proposition
shows.

\begin{proposition} Let $A$ be an $A_\infty$-algebra, $V$ a complex
and $f_1: A \to V$ a quasi-isomorphism of complexes. Then $V$ admits
a structure of $A_\infty$-algebra such that $f_1$ extends to
an $A_\infty$-quasi-isomorphism $f: A \to V$.
\end{proposition}

The analogous statement for dg algebras and their morphisms is
of course completely wrong. For our complex $V$, we can take in
particular the graded space $H^* V$ with the zero differential
(since we work over a field, the canonical surjection from the cycles
to the homology of $V$ splits and we obtain $f_1$ by composing
a right inverse with the inclusion of the cycles). Then we obtain
the first part of the

\begin{theorem} If $A$ is an $A_\infty$-algebra, then $H^*A$ admits
an $A_\infty$-algebra structure such that
\begin{itemize}
\item[(1)] $b_1=0$ and $b_2$ is induced from $b_2^A$, and
\item[(2)] there is an $A_\infty$-quasi-isomorphism $A \to H^*A$ inducing
the identity in homology.
\end{itemize}
Moreover, this structure is unique up to (non unique)
$A_\infty$-isomorphism.
\end{theorem}

Note that uniqueness up to $A_\infty$-quasi-isomorphism is trivial.
The point is that here we can omit `quasi'. An $A_\infty$-algebra
is {\em minimal} if $b_1=0$. The {\em minimal model} of an
$A_\infty$-algebra $A$ is the space $H^* A$ endowed with `the'
structure provided by the theorem. This structure can be
computed as follows \cite{KontsevichSoibelman00}: Choose
$$
\xymatrix{
A \ar@(ul,dl)[]_h \ar@<1ex>[r]^p & H^* A \ar@<1ex>[l]^i }
$$
such that $p$ and $i$ are morphisms of complexes of degree
$0$ and $h$ is a homogeneous map of degree $-1$ such that
$$
pi=\id \ko ip=\id + d(h) \ko h^2=0.
$$
Then the $n$th multiplication of the minimal model is
constructed as
$$
b_n^{min} = \sum_T {b_n}^T
$$
where $T$ ranges over the planar rooted trees $T$ with
$n$ leaves and $b_n^T$ is given by composing the tree-shaped
diagram obtained by labelling each leaf by $i$, each
branch point with $m$ branches by $b_m$, each internal
edge by $h$ and the root by $p$.
\medskip
$$
\centeps{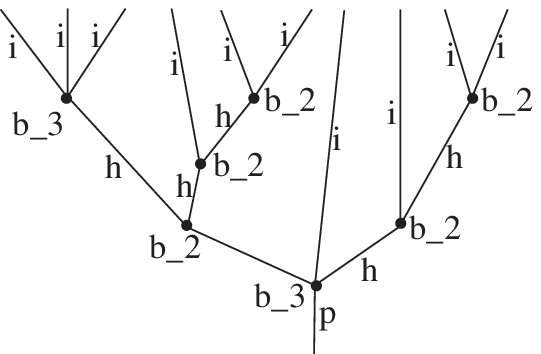}
$$

\subsection{Yoneda algebras} Let $B$ be a unital associative
algebra, $M$ a (right) $B$-module and $P \to M$ a projective resolution.
Let $A=\Hom_B(P,P)$ be the differential graded endomorphism
algebra of $P$ (its $n$th component consists of the morphisms
of graded objects of degree $n$ and its differential is the
supercommutator with the differential of $P$). Then $A$ is
in particular an $A_\infty$-algebra and thus has a minimal
model. Now the homology $H^* A$ is isomorphic, as an algebra
for $m_2$, to the Yoneda algebra
$$
\Ext^*_B(M,M).
$$
Thus we obtain higher multiplications on the Yoneda algebra.
The simplest case where these are non-trivial is that of
the algebra $B$ given by the quiver
$$
\xymatrix{
1 \ar[r]^\gamma & 2 \ar[r]^\beta & 3 \ar[r]^\alpha & 4
}
$$
with the relation $\alpha\beta\gamma=0$ and the module $M$
equal to the sum of the four simple $B$-modules. Then the
Yoneda algebra is given by the quiver
$$
\xymatrix{
1 & 2 \ar[l]_c & 3 \ar[l]_b & 4 \ar[l]_a \ar@(ul,ur)[lll]_e
} \ko
$$
where the arrows $a,b,c$ are of degree $1$, the arrow $e$
is of degree $2$, we have $m_2(c,b)=0$, $m_2(b,a)=0$
and $b_3(c,b,a)=e$.

Another example is given by the algebras $B=k[x]/(x^n)$, where
$n$ is an integer $\geq 3$, and the simple $B$-module $M=k$.
In this case, the underlying graded algebra of the Yoneda
algebra is isomorphic to $k[u,v]/(v^2)$, where $u$ is of
degree $2$ and $v$ of degree $1$. Note that this is
independent of $n$. In contrast, the $A_\infty$-structure
depends on $n$: it is given by the operations $b_l$ which
vanish for $l\neq 2,n$ such that $b_2$ is given
by the Yoneda multiplication $m_2$
and $b_n$ is the unique $k[u]$-multilinear map which vanishes
if one of the arguments is in $k[u]$ and satisfies
\[
b_n(v, \ldots, v)= u.
\]
This example was treated by D.~Madsen \cite{Madsen02}.
More generally, J.~Chuang and A.~King have computed the
$A_\infty$-structure on $Ext^*_B(M,M)$
for the sum $M$ of the simple modules over an
arbitrary monomial algebra $B$, \cf \cite{ChuangKingxy}.
Other examples have been computed by Lu-Palmieri-Wu-Zhang
\cite{LuPalmieriWuZhang04}.

\subsection{On the meaning of the higher
multiplications} \label{ss:meaning} This section should be
skipped on a first reading.
Suppose that $A$ is a minimal $A_\infty$-algebra. Endow $H^*(A)$
with the associative product obtained from $b_2$. Then $b_3$
yields a $3$-cocycle homogeneous of degree $-1$ for the natural
grading on $H^* A$ and thus an element in Hochschild cohomology
\[
\gamma \in HH^{3,-1}(H^* A, H^* A).
\]
This class determines the Massey products
in $H^*(A)$, \cf for example \cite{May69}. It also controls the
realizability (up to direct factors) of graded $H^*A$-modules
$M$ as the homology $H^* \tilde{M}$ of an $A$-module $\tilde{M}$,
\cf \cite{BensonKrauseSchwede04}.

One can also interpret the `higher' operations $b_n$. To simplify
the exposition, let us change the setting slightly: Let $\cb$
be a small $k$-linear abelian category and consider its
bounded derived category $\der^b(\cb)$. Thanks to the shift
functor, it becomes a graded category with morphism spaces
\[
\Homdot(L,M) = \bigoplus_{n\in\Z} \Hom_{\der^b(\cb)}(L, M[n]).
\]
This category can be completed into a minimal $A_\infty$-category
in much the same way that a Yoneda algebra can be completed into
a minimal $A_\infty$-algebra. In particular, the operation
$b_3$ determines an element in the Hochschild cohomology
\[
\gamma \in HH^{3,-1}(\der^b(\cb), \der^b(\cb))
\]
of the graded category $\der^b(\cb)$ (defined as in \cite{Mitchell72}).
One can show that this element determines the triangulated structure
of the $k$-linear category $\der^b(\cb)$. But in fact, it
determines more, namely it allows one to reconstruct the $k$-linear
category
\[
\der^b(\fun(\cdot\to\cdot, \cb))
\]
from $\der^b(\cb)$, where $\cdot\to\cdot$ denotes the $k$-category
with $2$ objects and one non trivial morphism and $\fun$ the
functor category. For this, one notes that $\gamma$ corresponds
to $\der^b(\fun(\cdot\to\cdot,\cb))$ considered as an extension of
$\fun(\cdot\to\cdot, \der^b(\cb))$ by the
bimodule (in the sense of \cite{Mitchell72})
\[
(L\to M) \mapsto \Hom_{\der^b(\cb)}(L, M[1]).
\]
In a similar manner, for $n\geq 1$, the operations $b_2, \ldots, b_{n+2}$
allow one to reconstruct
\[
\der^b(\fun( (\cdot\to\cdot)^n, \cb))
\]
from $\der^b(\cb)$.

\subsection{Units} \label{ss:Units}
A {\em strict unit} for an $A_\infty$-algebra $A$ is an element
$1\in A^0$ which is a unit for $m_2$ and such that, for $n\neq 2$,
the map $b_n$ takes the value $0$ as soon as one of its arguments
equals $1$. Unfortunately, strict unitality is not preserved
by $A_\infty$-quasi-isomorphisms. A {\em homological unit} for
$A$ is a unit for the associative algebra $H^*A$ with the multiplication
induced by $m_2$. Homological unitality is clearly preserved under
$A_\infty$-quasi-isomorphism but is not easy to handle in practical
computations. Fortunately, it turns out that the two notions are
not very different:

\begin{proposition}[\protect{\cite[3.2.1]{Lefevre03}}] Each (resp. minimal) homologically unital
$A_\infty$-algebra is $A_\infty$-quasi-isomorphic (resp. $A_\infty$-isomorphic)
to a strictly unital $A_\infty$-algebra.
\end{proposition}

\section{A-infinity modules} \label{s:A-infinity-modules}

Let $A$ be a homologically unital $A_\infty$-algebra. An
{\em $A_\infty$-module} is a graded space $M$ with maps
$$
b_n : SM \ten (SA)\tp{n-1} \to SM \ko n\geq 1 \ko
$$
homogeneous of degree $1$ such that the identity~\ref{ainfeq}
holds for all $n\geq 1$ (where we have to interpret $b_n$
as $b_n^A$ or $b_n^M$ according to the type of its arguments)
and that the induced action
$$
H^*M \ten H^*A \to H^*M
$$
is unital. For example, the $A_\infty$-algebra $A$ can
be viewed as a module over itself: the free module of
rank one. The notions of {\em morphism} and {\em quasi-isomorphism}
of $A_\infty$-modules are defined in the natural way. The
{\em derived category $\der_\infty A$} is defined as the localization
of the category of $A_\infty$-modules (with degree $0$ morphisms)
with respect to the class of quasi-isomorphisms. Thus, its objects
are all $A_\infty$-modules and its morphisms are obtained from
morphisms of $A_\infty$-modules by formally inverting all
quasi-isomorphisms. It turns out that the derived category
is naturally a triangulated category. The {\em perfect derived
category $\per A$} is defined as the closure of the free
$A$-module of rank one under shifts in both directions,
extensions and passage to direct factors.

When $A$ is an ordinary unital associative algebra, we
have a natural functor from the category of complexes
over the category $\Mod A$ of (right) $A$-modules to the
category of $A_\infty$-modules.
This functor is faithful but neither full nor essentially
surjective. Nevertheless, it induces an equivalence
$$
\der(\Mod A) \to \der_\infty A.
$$
Under this equivalence, the perfect derived category
corresponds to the full subcategory of complexes quasi-isomorphic
to bounded complexes of finitely generated projective
$A$-modules.

Let $B$ be another homologically unital $A_\infty$-algebra.
An {\em $A_\infty$-bimodule} is given by a graded space $X$
with maps
$$
b_{i,j} : (SA)\tp{i} \ten SX \ten (SB)\tp{j} \to SX \ko i+j\geq 0\ko
$$
satisfying the identity~\ref{ainfeq} for all $i,j\geq 0$
and such that $H^*X$ becomes a unital $H^*A$-$H^*B$-bimodule. For such
a bimodule, one can define the {\em tensor product}
$$
?\tensinf_A X : \der_\infty A \to \der_\infty B \ko
M \mapsto \bigoplus_{i=0}^\infty M \ten (SA)\tp{i} \ten X
$$
(we only indicate the underlying graded space).
The $A_\infty$-algebras $A$ and $B$ are {\em derived equivalent}
if there is an $A_\infty$-bimodule $X$ such that the associated
tensor product is an equivalence.

This generalizes the now classical notion of derived equivalence
for ordinary algebras: If $A$ and $B$ are ordinary algebras, then,
by J.~Rickard's theorem \cite{Rickard89}, they are derived equivalent iff $B$
admits a tilting complex (\eg a tilting module) with endomorphism
ring $A$.

Let us call {\em algebraic} a triangulated category which is the
homotopy category associated with a stable
$k$-linear Quillen model category, \cf section~\ref{ss:model-categories}
and \cite{Hovey99}
(recall that $k$ is the ground field).
The class of algebraic triangulated categories contains all
homotopy categories of complexes over $k$-linear categories and is stable
under passage to triangulated subcategories and Verdier localizations.
Thus, it contains all triangulated categories \lq of algebraic
origin\rq, \eg derived categories of categories of coherent sheaves.

Recall that an additive category has {\em split idempotents} if every
idempotent endomorphism admits a kernel (and thus gives rise to
a direct sum decomposition). By a {\em generator} of a triangulated
category $\ct$, we mean an object $G$ whose closure under shifts
in both directions, extensions and passage to direct factors
equals $\ct$.

\begin{theorem}[\protect{\cite[7.6]{Lefevre03}}] Let $\ct$ be a
($k$-linear) algebraic triangulated category
with split idempotents and a generator $G$. Then there is a structure
of $A_\infty$-algebra on
$$
A= \bigoplus_{n\in\Z} \Hom_\ct(G,G[n])
$$
such that $m_1=0$, $m_2$ is given by composition and that the
functor
$$
\ct \to \Grmod (A,m_2) \ko U \mapsto \bigoplus_{n\in\Z} \Hom_\ct(G, U[n])
$$
lifts to a triangle equivalence
$$
\ct \to \per(A).
$$
\end{theorem}

Here $\Grmod$ denotes the category of graded right modules.
For example, suppose that $\ct$ is the bounded derived category
of the category of coherent sheaves
on projective $n$-space. Then, as Beilinson has shown
\cite{Beilinson78} by \lq resolving
the diagonal\rq, the object $G=\bigoplus_{i=0}^n O(-i)$
generates $\ct$. Moreover, the algebra $A$ is concentrated in degree $0$
and of finite dimension and finite global dimension. Then we
obtain equivalences
$$
\der^b(\coh \bP^n) \iso \per(A) \liso \der^b(\mod A) \ko
$$
where $\mod A$ denotes the category of finite-dimensional $A$-modules.

Beautiful theorems on the existence of generators in triangulated
categories of geometric origin are due to Bondal-Van den Bergh
\cite{BondalVandenBergh03} and to Rouquier \cite{Rouquier03}.

\section{The bar-cobar equivalence}

In this section, we present the interpretation of
$A_\infty$-algebras as the fibrant objects in the model category of
certain differential graded coalgebras following \cite{Lefevre03}.
For the reader's convenience we first recall the notion of
model category.

\subsection{Reminder on model categories}
 \label{ss:model-categories}
Model categories were invented by D.~Quillen \cite{Quillen67}.
A good modern introduction can be found in \cite{DwyerSpalinski95}.
A {\em model category} is a category $\cc$ endowed with three
classes of morphisms
\begin{itemize}
\item[-] the class $\cw$ whose elements are called {\em weak equivalences}
and denoted by arrows with tilde $\iso$,
\item[-] the class $\fib$ whose elements are called {\em fibrations} and
denoted by feathered arrows $\xymatrix{ \ar@{>->}[r] & }$,
\item[-] the class $\cof$ whose elements are called {\em cofibrations}
and denoted by double-headed arrows $\xymatrix{ \ar@{->>}[r] & }$.
\end{itemize}
Each of these classes is supposed to be stable under compositions
and to contain the identities. Moreover, the classes are supposed
to satisfy the following axioms
\begin{itemize}
\item[MC1] The category $\cc$ admits all finite limits and colimits
(in particular, it admits an initial object $\emptyset$ and a terminal
object $*$).
\item[MC2] The class $\cw$ is {\em saturated}, \ie if two among
$f$, $g$, and $fg$ belong to $\cw$ so does the third.
\item[MC3] The three classes $\cw$, $\fib$ and $\cof$ are stable
under retracts in the category of morphisms of $\cc$.
\item[MC4] In each commutative square of solid arrows of the
following two types
\[
\newdir{ >}{{}*!/-7pt/\dir{>}}
\xymatrix{  \ar[r] \ar@{ >->}[d]_{i}^{\tilde{}} & \ar@{->>}[d]^p \\
\ar[r] \ar@{.>}[ru] & }
\quad\quad\quad
\xymatrix{  \ar[r] \ar@{ >->}[d]_j & \ar@{->>}[d]_{\tilde{ }}^q \\
\ar[r] \ar@{.>}[ru] & }
\]
the dotted arrow exists so as to make the two triangles commutative.
\item[MC5] Each morphism $f$ factors as $f=pi$, where $p$ is a fibration
and $i$ a cofibration and a weak equivalence. It also factors as
$f=qj$, where $q$ is a fibration and a weak equivalence and $j$
a cofibration.
\end{itemize}
Note that these axioms are selfdual. In the situation of MC4, one says
that $p$ has the {\em right lifting property} with respect to $i$ and
$i$ the {\em left lifting propery} with respect to $p$. It is not hard
to see that if $\cc$ is a model category, then the class $\fib$ is
determined by the class of {\em trivial cofibrations} $\cw \cap \cof$
and the class $\cof$ is determined by the class of {\em trivial
fibrations $\cw \cap\fib$}.  In fact, $\fib$ consists of the morphisms
with the right lifting property for all trivial cofibrations and
$\cof$ of the morphisms with the left lifting property with respect to
all trivial fibrations.  A model category is {\em left proper} if any
pushout of a weak equivalence along a cofibration is a weak
equivalence.

As an example, consider the category $\cc=\cc^+(\Mod B)$ of left
bounded complexes
\[
\ldots \to 0 \to \ldots \to X^p \to X^{p+1} \to \ldots
\]
of right modules over a ring $B$. We define $\cw$ to be the
class of quasi-isomorphisms, $\cof$ to be the class of
morphisms $i: X \to Y$ such that each $i^n$, $n\in\Z$, is injective
and $\fib$ to be the class of morphisms $p: X\to Y$ such
that each morphism $p^n$, $n\in\Z$, is surjective and has
an injective kernel. It is not hard to show, \cf \cite{Hovey99},
that $\cc$ is a model category. Note that in $\cc$ the
initial and the terminal object coincide with the complex
of zero modules, that the morphism $0 \to X$ is always a
cofibration and that the morphism $X \to 0$ is a fibration
iff all components $X^n$, $n\in\Z$, are injective. Therefore,
in this example, the first half of axiom MC5 corresponds
to the existence of injective resolutions for complexes
\[
\newdir{ >}{{}*!/-7pt/\dir{>}}
\xymatrix{ X \ar[rr] \ar@{ >->}[rd]_{\tilde{ }} &  & 0 .\\
 & I \ar@{->>}[ur]}
\]

In an arbitrary model category $\cc$, an object $X$ is
{\em fibrant} if the morphism $X \to *$ is a fibration;
an object $Y$ is {\em cofibrant} if the morphism $\emptyset \to Y$
is a cofibration.

Thus, in the above example, all complexes $X$ are cofibrant and a
complex $Y$ is fibrant iff it has injective components.

For an arbitrary model category $\cc$, the {\em homotopy category}
$\Ho(\cc)$ is defined as the localization $\cc[\cw^{-1}]$ of
the category $\cc$ with respect to the class $\cw$. Thus, it
has the same objects as $\cc$ and its morphisms are obtained
from the morphisms of $\cc$ by formally inverting all elements
of $\cw$. In the above example, the homotopy category equals
the left bounded derived category $\der^+(\Mod B)$.

\begin{theorem} Let $\cc$ be a model category and $\cc_{cf}$ the full
subcategory formed by the objects $X$ which are both fibrant and
cofibrant. Then, for $X,Y\in\cc_{cf}$, the canonical map
\[
\Hom_{\cc} (X,Y) \to \Hom_{\Ho(\cc)}(X,Y)
\]
is surjective.
\end{theorem}
The proof of the theorem can be found in \cite{DwyerSpalinski95}, for
example.
With the notations of the theorem, define two morphisms $f,g$
from $X$ to $Y$ to be {\em homotopic} if their images in the
homotopy category coincide. Is is then clear that the canonical
functor
\[
\cc_{cf}/\mbox{homotopy} \to \Ho(\cc)
\]
is an equivalence.

In our example, a complex $X$ belongs to $\cc_{cf}$ iff it
has injective components. Two morphisms $f,g$ between such complexes
are then homotopic in the sense of the model category iff they
are homotopic in the classical sense, \ie there is a
$B$-linear morphism $r: X \to Y$ homogeneous of degree $-1$ such
that
\[
f-g = d \circ r + r \circ d.
\]
Finally, the above equivalence specializes to the classical
equivalence
\[
\cc^+(\Inj B)/\mbox{homotopy} \iso \der^+(\Mod B) \ko
\]
where $\Inj B$ denotes the full subcategory of injective
$B$-modules.

In a model category $\cc$, an object $I$ is {\em minimal}
if it is fibrant and cofibrant and for all endomorphisms
$f$ of $I$ in $\cc$, we have
\[
f \mbox{ is invertible in } \cc \Leftrightarrow
f \mbox{ becomes invertible in } \Ho(\cc).
\]
A {\em minimal model} for $X\in \cc$ is an isomorphism
\[
X \iso I
\]
of $\Ho(\cc)$ such that $I$ is minimal.

In our example, assume that $X=M$ is a $B$-module, \ie
a complex concentrated in degree $0$. Suppose that we have
an injective resolution
\[
0 \to M \to I^0 \to I^1 \to \ldots .
\]
Such a resolution gives rise to an isomorphism $M \iso I$
in the derived category, \ie $\Ho(\cc)$. This is a minimal
model for $M$ iff the above resolution is a minimal injective
resolution.

\subsection{Algebras}
Let $k$ be a field. Let $\alg$ be the category of augmented
dg algebras, \ie unital associative dg algebras $A$ endowed with a
morphism
\[
\eps : A \to k
\]
of unital dg algebras. In $\alg$, define $\cw$ to be the class of
quasi-isomorphisms and $\fib$ to be the class of surjective
morphisms. Let $\cof$ be the class of morphisms $i$ such that
in each commutative square of solid arrows in $\alg$
\[
\newdir{ >}{{}*!/-7pt/\dir{>}}
\xymatrix{  \ar[r] \ar@{ >->}[d]_i & \ar@{->>}[d]^{\tilde{ }} \\
\ar[r] \ar@{.>}[ru] & }
\]
there is a dotted arrow which makes the two triangles
commutative. Note that the algebra $k$ is both the initial
and the terminal object of $\alg$.

\begin{theorem}[Hinich \cite{Hinich97}] The category $\alg$ endowed
with the above classes $\cw$, $\fib$ and $\cof$ is a model category.
All its objects are fibrant.
\end{theorem}

Below (in Theorem~\ref{coalgebra-model}),
we will give a description of the cofibrant objects of
$\alg$. We will see that if an object of $\alg$ is cofibrant, then its
underlying graded algebra is free.
However, one can show that
this condition is not sufficient. One can also show that not all
objects $A\in \alg$ admit minimal models.

\subsection{Coalgebras and the bar/cobar constructions}
Let $\cog'$ be the category of augmented dg coalgebras.
Thus, an object of $\cog'$ is a complex $C$ endowed with a
coassociative comultiplication
\[
\Delta : C \to  C \ten C
\]
for which $d$ is a {\em coderivation}, \ie
\[
\Delta\circ d = (d\ten \id + \id \ten d) \circ \Delta \ko
\]
such that $\Delta$ admits a counit $\eta : C \to k$ and a (co-)augmentation
$\eps : k \to C$, \ie a counital coalgebra morphism such that $\eta\eps =
\id_k$. For example, if $V$ is a complex, then we let $T^c(V)$ be the
tensor coalgebra on $V$: Its underlying complex is
\[
\bigoplus_{n\geq 0} V\tp{n}
\]
and its comultiplication `separates tensors':
\[
\Delta(v_1, \ldots, v_n) = \sum_{i=0}^n (v_1, \ldots, v_i) \ten
(v_{i+1}, \ldots, v_{n}) \ko
\]
where $(v_1, \ldots, v_n)$ stands for $v_1\ten \ldots\ten v_n$
and the empty parentheses $()$ are to be interpreted as $1$.
Then $T^c(V)$ is endowed with a canonical projection
\[
T^c(V) \to V
\]
dual to the inclusion $V \to T(V)$ of $V$ into the tensor algebra,
but in the category $\cog'$, the coalgebra $T^c(V)$ does not have
the universal property dual to that of the tensor algebra: it is
not cofree on $V$. To resolve this technical problem is one
of the motivations for replacing $\cog'$ by a smaller category,
namely that of cocomplete augmented dg coalgebras. A coalgebra
$C\in \cog'$ is {\em cocomplete} if it is the union of the
kernels of the maps
\[
C \to C\tp{n} \to (C/k)\tp{n} \ko n\geq 2\ko
\]
which are the compositions of the canonical projection with the
iterated comultiplications. Define $\cog$ to be the full subcategory
of $\cog'$ whose objects are the cocomplete coalgebras $C\in \cog'$.
Then it is not hard to show that $T^c(V)$ is cocomplete
and that it is cofree on $V$ in the category $\cog$.

For a dg coalgebra $C$ and a dg algebra $A$, consider the complex
\[
\Homdot_k(C,A)
\]
whose $n$th component is the space of homogeneous $k$-linear
maps $f:C \to A$ of degree $n$ and whose differential maps $f$ to
\[
d\circ f - (-1)^n f\circ d.
\]
This complex becomes a dg algebra for the {\em convolution} defined by
\[
f*g = \mu \circ (f\ten g) \circ \Delta.
\]
Define a homogeneous $k$-linear map $\tau : C \to A$ to be a
{\em twisting cochain} if it is homogeneous of degree $1$ and satisfies
\[
d(\tau) + \tau*\tau = 0 \mbox{ and } \eps\circ\tau\circ\eps=0.
\]
Let $\Tw(C,A)$ be the set of twisting cochains. Then it is not
hard to show that for a given $A\in \alg$, the functor
\[
\cog \to \sets \ko C \mapsto \Tw(C,A)
\]
is representable. The {\em bar construction $BA$} is by definition
a representative. It is endowed with the {\em universal twisting cochain
$\tau_0$} which corresponds to the identity of $BA$ under the bijection
\[
\Tw(BA,A) = \Hom_{\cog}(BA,BA).
\]
One checks that we have $BA=T^c(SA)$ endowed with
a differential which takes into account the differential
of $A$ and its multiplication. The universal twisting cochain
is the canonical projection $BA \to A$.

Dually, one shows that for a coalgebra $C$, the functor
\[
\alg \to \sets \ko A \mapsto \Tw(C,A)
\]
is corepresentable. The {\em cobar construction $\Omega(C)$} is
a corepresentative. One has $\Omega C= T(S^{-1}C)$ and the
differential of $\Omega C$ takes into account the differential
of $C$ and the comultiplication. The universal twisting cochain is
the canonical inclusion $C \to \Omega C$.

It is clear from these constructions that for $A\in\alg$ and $C\in \cog$,
we have
\[
\Hom_{\alg}(\Omega(C), A) = \Tw(C,A) = \Hom_{\cog}(C, B(A))
\]
so that $\Omega$ is left adjoint to $B$. If $A$ is an ordinary augmented
algebra (concentrated in degree $0$) and $C$ an ordinary augmented
coalgebra, one checks that
\[
H_*(BA) = \Tor_*^A(k,k) \mbox{ and } H^*(\Omega C) = \Cotor^*_C(k,k).
\]
In $\cog$, we define $\cw$ to be the class of morphisms $f: C \to
D$ such that $\Omega(f)$ is a quasi-isomorphism; we define $\cof$
to be the class of injective morphisms and $\fib$ to be the class
of morphisms $p$ such that in each commutative square of solid
arrows in $\cog$
\[
\newdir{ >}{{}*!/-7pt/\dir{>}}
\xymatrix{  \ar[r] \ar@{ >->}[d]_{\tilde{ }} & \ar@{->>}[d]^{p} \\
\ar[r] \ar@{.>}[ru] & }
\]
there is a dotted arrow which makes the two triangles commutative.

\begin{theorem}[\protect{Lef\`evre \cite[1.3]{Lefevre03}}]
\label{coalgebra-model}
\begin{itemize}
\item[a)] The category $\cog$ endowed with the
above classes $\cw$, $\cof$ and $\fib$ is a model category.
The functors $B$ and $\Omega$ induce quasi-inverse equivalences
\[
\Ho(\alg) \iso \Ho(\cog).
\]
Moreover, $B$ preserves fibrations and trivial fibrations and
$\Omega$ preserves cofibrations and trivial cofibrations.
\item[b)] All objects of $\alg$ are fibrant. An object of $\alg$ is cofibrant
iff it is a retract of $\Omega C$ for some $C$ in $\cog$.
All objects of $\cog$ are cofibrant. An object is fibrant
iff its underlying graded coalgebra is a tensor coalgebra.
\item[c)] Two morphisms $f,g: A \to A'$ between
fibrant-cofibrant objects of $\alg$ are homotopic iff there is a $k$-linear
map $h: A \to A'$ homogeneous of degree $-1$ such that
\[
\eps\circ h =0 \mbox{ and }
h\circ m_A = m_{A'} \circ (f \ten h + h \ten g) \mbox{ and }
f-g = d\circ h + h \circ d.
\]
The dual statement holds for morphisms between fibrant-cofibrant
objects of $\cog$.
\end{itemize}
\end{theorem}

Notice that $\alg$ and $\cog$ behave similarly in many respects
but that the descriptions of their fibrant-cofibrant objects are
quite different. Each cofibrant dg algebra has a free underlying
graded algebra (since retracts of free algebras are free).
Below, in
section~\ref{cofibrant-algebras}, we give an example
of a non cofibrant object of $\alg$ whose underlying graded
algebra is free. Hinich \cite{Hinich97} describes the
cofibrations of $\alg$ as `transfinite cell attachments'.

\subsection{The link with $A_\infty$-algebras and minimal models}
An {\em augmented} $A_\infty$-algebra is a strictly unital
$A_\infty$-algebra endowed with a strict unit-preserving
morphism $\eps: A \to k$.
Such an algebra decomposes as $A=k\oplus \ol{A}$ and the
functor $A \mapsto \ol{A}$ yields an equivalence between the
category of augmented $A_\infty$-algebras and the category
of non unital $A_\infty$-algebras.

For a coalgebra $C\in \cog$ and an augmented
$A_\infty$-algebra $A$, the complex $\Homdot_k(C,A)$ becomes an
augmented $A_\infty$-algebra for the convolution operations
\[
b_n (f_1, \ldots, f_n) = b_n^A \circ (f_1 \ten \cdots\ten f_n)
\circ \Delta^{(n)}\ko
\]
where $\Delta^{(n)}$ is the iterate of $\Delta$ taking values in
$C\tp{n}$. Define $\Tw_\infty(C,A)$ to be the set of solutions
$\tau$ of the `Maurer-Cartan equation'
\[
\sum_{n\geq 1} b_n(\tau,\ldots,\tau) =0 \ko
\]
\cf \cite{Keller01}.
Note that the sum yields a well-defined map $C \to A$ thanks to
the fact that $C$ is cocomplete. Then the functor
\[
\cog \to \sets \ko C \to \Tw_\infty(C,A)
\]
is representable and we denote by $B_\infty A$ a representative.
It is not hard to check that $B_\infty A$ is $T^c(SA)$ endowed
with the unique coderivation whose composition with the projection
$BA \to SA$ has the components
\[
b_n: (SA)\tp{n} \to SA \ko n \geq 1.
\]
In fact, if $V$ is a graded vector space, then there is a
bijection between $A_\infty$-structures on $V$ and coalgebra
differentials on $TSV$. Moreover, if $A$ and $A'$ are
$A_\infty$-algebras, then the $A_\infty$-morphisms
$f: A \to A'$ are in canonical bijection with the
morphisms $B_\infty A \to B_\infty A'$ of augmented
dg coalgebras.

\begin{theorem}[\protect{\cite[1.3]{Lefevre03}}]
\begin{itemize}
\item[a)] An object of $\cog$ is fibrant-cofibrant iff it is
isomorphic to $B_\infty A$ for some $A_\infty$-algebra $A$.
\item[b)] Each object $C$ of $\cog$ admits a minimal model and, if
$A$ is an $A_\infty$-algebra and $A_{min}$ its minimal model, then
$B_\infty(A_{min})$ is the minimal model of $B_\infty A$.
\end{itemize}
\end{theorem}

It is not hard to see that an $A_\infty$-quasi-isomorphism $f:A \to
A'$ between $A_\infty$-algebras yields a weak equivalence
$B_\infty A \to B_\infty A'$. Since both of these coalgebras are
fibrant-cofibrant, such a weak equivalence has an inverse up to
homotopy. This yields that $f$ has an inverse up to
$A_\infty$-homotopy. The theorem also yields a good interpretation
of the notion of minimality for $A_\infty$-algebras. The following
square summarizes the situation. We denote the category of
augmented $A_\infty$-algebras by $\alg_\infty$.
\[
\xymatrix{
\alg \ar@<3pt>[d]^B \ar@{^{(}->}[rr]_{\mbox{{\scriptsize not full}}} &  &\alg_\infty \ar[d]^{B_\infty}_{\tilde{ }} \\
\cog \ar@<3pt>[u]^\Omega & & \cog_{cf} \ar[ll]^{\mbox{{\scriptsize fully faithful}}}
}
\]
If we pass to the homotopy categories (where homotopy means
localization for $\alg$ and $\cog$ and passage to quotient
categories for $\alg_\infty$ and $\cog_{cf}$), then all the
arrows of the diagram become equivalences. Thus we obtain four
descriptions of the localization of the category of dg algebras
with respect to all quasi-isomorphisms. In particular,
$A_\infty$-algebras do not yield `new homotopy types' of algebras
but a new description of the existing types.

\subsection{A non cofibrant dg algebra}
\label{cofibrant-algebras}
Here is an example of a dg algebra $A\in \alg$
which is free as a graded algebra but not cofibrant: Take $A=TV$,
where $V=k$ is concentrated in degree $1$. Endow $A$ with the
unique differential whose restriction to $V\subset TV$ is
\[
V = k \iso k\ten k = V\tp{2}\subset TV.
\]
Then clearly $A$ is quasi-isomorphic to its subalgebra $k$, which
is fibrant-cofibrant. If $A$ was fibrant-cofibrant as well, then
the inclusion $k\to TV$ should admit a left inverse up to
homotopy in the sense of c) in theorem~\ref{coalgebra-model}.
But there are no non-zero maps $h: A \to k$ of degree $-1$ such
that $\eps\circ h =0$. Therefore, $A$ cannot be
fibrant-cofibrant. Note that $A$ is the cobar construction
on a non cocomplete dg coalgebra.

\subsection{Modules and comodules} Let $A$ be an augmented
dg algebra and $C$ a cocomplete augmented dg coalgebra. Let
$\tau: C \to A$ be a twisting cochain. For a dg right $A$-module
$L$, we endow the $C$-comodule $L\ten C$ with the differential
defined by
\[
d = d_L\ten \id_C + \id_L\ten d_C +
(\mu\ten\id_C)(\id_L\ten\tau\ten\id_C)(\id_L\ten\Delta)\ko
\]
where $\Delta: C \to C\ten C$ is the comultiplication of $C$
and $\mu: L \ten A \to L$ the multiplication of $L$.
We write $L\ten_\tau C$ for the resulting dg $C$-comodule.
Similarly, if $M$ is a dg right $C$-comodule, we define
a differential on $M\ten A$ by
\[
d=d_M\ten \id_A + \id_M \ten d_A +
(\id_M\ten\mu_A)(\id_M\ten\tau\ten\id_A)(\delta_M\ten\id_A).
\]
We denote the resulting dg $A$-module by $M\ten_\tau A$.

Let us denote by $\Mod A$ the category of right dg $A$-modules and
by $\Comc C$ the category of dg $C$-comodules $M$ which are {\em cocomplete},
\ie $M$ is the union of the kernels of the maps induced by
iterated comultiplications
\[
M \to M \ten \ol{C}\tp{n} \ko n\geq 2.
\]
Then the pair
\[
\xymatrix{ \Mod A \ar@<3pt>[d]^{?\ten_\tau C} \\
\Comc C \ar@<3pt>[u]^{?\ten_\tau A}
}
\]
is a pair of adjoint functors (\cf Lemme 2.2.1.2 of \cite{Lefevre03}).
From now on, we suppose that $\tau$ is {\em acyclic}, \ie that the
following equivalent (\cf Proposition 2.2.4.1 of \cite{Lefevre03})
conditions hold
\begin{itemize}
\item[(i)] The adjunction morphism $M\ten_\tau C \ten_\tau A \to M$
is a quasi-isomorphism for each dg $A$-module $M$.
\item[(ii)] The adjunction morphism $A\ten_\tau C \ten_\tau A \to A$
is a quasi-isomorphism.
\item[(iii)] The morphism $BC \to A$ induced by $\tau$ is a quasi-isomorphism.
\item[(iv)] The morphism $C \to \Omega A$ induced by $\tau$ is a
weak equivalence.
\end{itemize}

If follows from \cite{Hinich97} that the category $\Mod A$ admits
a unique structure of model category whose weak equivalences are
the quasi-isomorphisms and whose fibrations are the surjective
morphisms. We write $\der(A)$ for the localization of $\Mod A$
with respect to the class of quasi-isomorphisms.

\begin{theorem}[Th\'eor\`eme 2.2.2.2 of \cite{Lefevre03}] \label{thm-modules}
\begin{itemize}
\item[a)] The category $\Comc C$ admits a unique structure of
model category whose weak equivalences are the morphisms
$f$ such that $f\ten_\tau A$ is a quasi-isomorphism and whose
cofibrations are the injetive morphisms. We write $\der(C)$
for the localization of $\Comc C$ with respect to the class
of weak equivalences.
\item[b)] The functors $?\ten_\tau C$ and $?\ten_\tau A$ induce
quasi-inverse equivalences
\[
\der(A) \iso \der(C).
\]
\end{itemize}
\end{theorem}

\subsection{Example: Koszul algebras} Let $V$ be a vector
space and $R\subset V\ten V$ a subspace such that the
algebra $A=TV/(R)$ is Koszul \cite{BeilinsonGinzburgSoergel96},
\ie the $A$-module $k$
admits a minimal projective resolution $P \to k$ in
the category of graded $A$-modules such that $P_i$
is generated in degree $i$ for all $i\geq 0$. Let
$C$ be the subspace of the tensor coalgebra $T^c(V)$
whose $n$th component is $k$ for $n=0$, $V$ for $n=1$ and
\[
\bigcap_{p+2+q=n} V\tp{p}\ten R \ten V\tp{q}
\]
for $n\geq 2$. Then clearly $C$ is a subcoalgebra.
We endow $C$ with the grading such that $V$ is in degree
$-1$ and with the vanishing differential. We consider
$A$ as a dg algebra concentrated in degree $0$. We let
$\tau: A \to C$ be given by the composition
\[
C \to V \to A
\]
of the inclusion with the projection. Clearly, $\tau$
is a twisting cochain. Then the complex
\[
A \ten_\tau C \ten_\tau A
\]
is the Koszul resolution of the bimodule $A$. Thus, the
twisting cochain $\tau$ is acyclic and we obtain an equivalence
of categories
\[
?\ten_\tau C: \der(A) \to \der(C)
\]
which takes $A$ to $k$ (weakly equivalent to $A\ten_\tau C$)
and $k$ to $C$. An analogous equivalence can be constructed between
suitable (Adams-) graded versions of the categories $\der(A)$
and $\der(C)$. We refer to \cite{Floystad00} for related developments.

\subsection{Link with $A_\infty$-modules} \label{link-Ainf-modules}
Let $A$ be an
augmented $A_\infty$-algebra. Let $C=B_\infty A$. The
adjunction morphism
\[
B_\infty A = C \to B \Omega C = B_\infty(\Omega C)
\]
corresponds to a canonical $A_\infty$-morphism $A \to \Omega B_\infty A$,
which is a quasi-isomorphism and which is universal among
the $A_\infty$-morphisms from $A$ to a dg algebra
(Lemme 2.3.4.3 of \cite{Lefevre03}). Let us
put $U(A)=\Omega B_\infty A$. Thus we have a canonical
acyclic twisting cochain $\tau: B_\infty(A) \to U(A)$.
By theorem~\ref{thm-modules}, we have an equivalence
\[
\der(U(A)) \iso \der(B_\infty(A)).
\]
Both of these categories can be regarded as the (unbounded)
derived category of the $A_\infty$-algebra $A$. The link
with $A_\infty$-modules is the following:
Let $\Mod_\infty A$ be the category of $A_\infty$-modules
over $A$ which are strictly unital (\ie $b^M_n$, $n\geq 2$,
vanishes as soon as one of the arguments is $1$) and whose
morphisms are strictly unital (\ie $f_n$, $n\geq 2$, vanishes
as soon as one of the arguments is $1$). This category is
isomorphic to the category of {\em all} $A_\infty$-modules
with {\em all} $A_\infty$-morphisms over the reduction $\ol{A}$.
By restriction along the $A_\infty$-morphism $A \to U(A)$,
each dg module over $U(A)$ yields an $A_\infty$-module
in $\Mod_\infty A$. However, the resulting functor
\[
\Mod U(A) \to \Mod_\infty A
\]
is not full (its image only contains strict morphisms).

Let $M$ be in $\Mod_\infty A$. The datum of its strictly
unital $A_\infty$-module structure over $A$ is equivalent
to the datum of an arbitrary $A_\infty$-module structure
over $\ol{A}$. This in turn is equivalent to giving a
comodule differential on the induced comodule $M\ten B_\infty A$.
We write $B_\infty M$ for the induced comodule endowed
with the differential corresponding to a given
$A_\infty$-module structure on $M$.
Thus we obtain a functor
\[
\Mod_\infty A \to \Comc B_\infty(A) \ko M \to B_\infty M.
\]
It is not hard to see that this functor is in fact
fully faithful. We define two morphisms of $\Mod_\infty A$
to be homotopic if they become homotopic in $\Comc B_\infty(A)$.

\begin{proposition}
\begin{itemize}
\item[a)] The functor $M \mapsto B_\infty M$
induces an equivalence onto the full subcategory
$(\Comc B_\infty A)_{cf}$ of fibrant-cofibrant objects
of $\Comc B_\infty(A)$.
\item[b)] The functor $M \mapsto B_\infty M$ induces an
equivalence
\[
(\Mod_\infty A)/\mbox{homotopy} \iso \der(C).
\]
\end{itemize}
\end{proposition}

\begin{proof} a) In Proposition 2.4.1.3 of \cite{Lefevre03}, it
is shown that $B_\infty M$ is fibrant (and cofibrant) in
$\Comc B_\infty A$ and that each fibrant object is a direct
factor of an object $M\ten_\tau B_\infty A$ for some $U(A)$-module
$M$. It remains to be shown that such a direct factor is
of the form $B_\infty M'$ for some $M'\in \Mod_\infty A$.
Put $C=B_\infty A$. We will use that a morphism $g: V\ten C \to W\ten C$
of induced $C$-comodules is invertible iff it induces an
invertible morphism in the spaces of primitive elements
\[
V=(V\ten C)_{[1]} = \ker(
\xymatrix{V\ten C \ar[rrr]^-{\id_V\ten \Delta -
\id_V\ten\id_V\ten\eta} & & & V\ten C \ten C}).
\]
If $f$ is an idempotent endomorphism of
$M\ten_\tau C$, it induces an idempotent endomorphism of
graded vector spaces in $M = (M\ten_\tau C)_{[1]}$. Let $M=M'\oplus M''$
be the corresponding decomposition in the category of graded
vector spaces. Then we have $M\ten C = (M'\ten C) \oplus (M'' \ten C)$
as $C$-comodules and in the corresponding decomposition
\[
f=\left[ \begin{array}{cc} f_{11} & f_{12} \\ f_{21} & f_{22} \end{array}
\right]
\]
the diagonal entries are invertible. It follows that the
image $I$ of $f$ is isomorphic to $M'\ten C$ as a $C$-comodule.
The differential of $M\ten_\tau C$ yields a comodule differential
on $I\iso M'\ten C$ and thus an $A_\infty$-module structure on
$M'$. We have $I\iso B_\infty M'$.
Point b) is proved in Lemme 2.4.2.3
of \cite{Lefevre03}.
\end{proof}

To sum up, we have the following diagram
\[
\xymatrix{
\Mod U(A) \ar@{^{(}->}[rr]^{\mbox{{\scriptsize not full}}}
\ar@<3pt>[d]^{?\ten_\tau B_\infty A} & &
\Mod_\infty A \ar[d]^{B_\infty}_{\tilde{ }} \\
\Comc B_\infty(A) \ar@<3pt>[u]^{?\ten_\tau U(A)} & &
(\Comc B_\infty A)_{cf} \ar[ll]_{\mbox{{\scriptsize fully faithful}}} }
\]
All functors induce equivalences in the localizations
(for the two categories on the left) respectively the quotients
by the homotopy relations (for the two categories on the right).

\subsection{The derived category of a non augmented $A_\infty$-algebra}
Let $A$ be a (non augmented) $A_\infty$-algebra. Let $A^+=A\oplus k$
be the augmented $A_\infty$-algebra obtained by adjoining $k$.
The augmentation $A^+ \to k$ yields a functor
\[
\Mod_\infty A^+ \to \Mod_\infty k
\]
which passes to the derived categories, \cf
section 4.1.1 of \cite{Lefevre03}. By definition, the
(unbounded) derived category $\der_\infty(A)$ is the kernel of the
functor
\[
\der_\infty(A^+) \to \der_\infty (k).
\]
One can show that a strictly unital $A_\infty$-module $M$
over $A^+$ is in the kernel iff $B_\infty M$ is acyclic
(Remarque 4.1.3.5 of \cite{Lefevre03}).

Now suppose that $A$ is
homologically unital, \ie $H^*(A)$ endowed with the
multiplication induced by $m_2$ is unital.  Then $M$
belongs to the kernel iff $M$ is homologically
unital, \ie $H^*(M)$ is a unital $H^*(A)$-module,
\cf Lemme 4.1.3.7 of \cite{Lefevre03}. Moreover,
in this case $\der_\infty(A)$ is a compactly generated
triangulated category and has the free $A$-module
of rank one as a compact generator. Indeed, $A^+$
is quasi-isomorphic to $U(A^+)$ and one checks
that if $B=k\oplus \ol{B}$ is an augmented dg algebra with
homologically unital augmentation ideal $\ol{B}$,
then $\ol{B}$ is a compact generator for the
kernel of $\der(B) \to \der(k)$.

\section{A conceptual construction of $A_\infty$-functor
categories}

The notion of $A_\infty$-algebra naturally generalizes to
that of $A_\infty$-category. It has been known for some time,
\cf \cite{Kontsevich98},
that for two $A_\infty$-categories $\ca$, $\cb$, there is
a natural $A_\infty$-category $\Funinf(\ca,\cb)$ and for three
$A_\infty$-categories, a natural  `$A_\infty$-bifunctor'
\[
\Funinf(\cb,\cc)\times\Funinf(\ca,\cb) \to \Funinf(\ca,\cc).
\]
One can use this composition bifunctor to give a natural
construction for the bialgebra structure of \cite{Baues81},
\cite{GetzlerJones94} on the bar
construction $BC$ of the Hochschild complex $C=C(A,A)$
of an algebra $A$, \cf section~\ref{Ainf-functors}.

The $A_\infty$-category $\Funinf(\ca,\cb)$ was constructed in
\cite{Lefevre03} using twists of $A_\infty$-structures.  Here, we
further develop an idea of V.~Lyubashenko's \cite{Lyubashenko03} and
interpret $\Funinf(\ca,\cb)$ as an internal $\Hom$-object in a tensor
category, namely the category of cocomplete augmented dg
cocategories. In their work \cite{BespalovLyubashenkoManzyuk06},
Yu.~Bespalov, V.~Lyubashenko and O.~Manzyuk are developing a
comprehensive theory based on this interpretation. The main difference with
V.~Lyubashenko's original approach in \cite{Lyubashenko03} is that he
only considered the tensor subcategory generated by all cofree
cocategories. The construction of $\Funinf$ given in 
\cite{BespalovLyubashenkoManzyuk06} agrees with the one
below.

\subsection{Augmented $k$-quivers} \label{ss:Augmented-k-quivers}
Let $k$ be a commutative ring. A {\em $k$-quiver $V$} consists of
a set of {\em objects} $\obj(V)$ and of $k$-modules $V(x,y)$, for
all objects $x,y$ of $V$. We call the elements of $V(x,y)$ the
{\em morphisms} from $x$ to $y$ in $V$. A {\em morphism of
$k$-quivers} $F: V \ra W$ consists of a map $F: \obj(V)\ra\obj(W)$
and of maps
\[
F_{x,y} : V(x,y) \ra W(Fx,Fy)\ko
\]
for all objects $x,y$ of $V$. For each set $S$, we have the {\em
discrete $k$-quiver $kS$} with $\obj(kS)=S$, $(kS)(x,x)=k\id_x$ and
$(kS)(x,y)=0$ for all objects $x\neq y$ of $kS$. An {\em augmented
$k$-quiver} is a $k$-quiver $V$ endowed with morphisms
\[
k \obj(V) \arr{\eta} V \arr{\eps} k\obj(V)
\]
whose composition is the identity of $k\obj(V)$. A {\em morphism of
augmented $k$-quivers} is a morphism of the underlying quivers compatible
with the morphisms $\eta$ and $\eps$. The {\em tensor product} of
two augmented $k$-quivers $V$ and $W$ is the naturally augmented
$k$-quiver $V\ten W$ whose
set of objects is $\obj(V)\times \obj(W)$ and such that
\[
(V\ten W)\,((x,x'),(y,y'))= V(x,y)\ten W(x',y').
\]
The {\em unit} augmented $k$-quiver $E$ has one object $*$ and
$E(*,*)=k$. It is endowed with the identity morphisms $\eta$ and
$\eps$. Clearly, we have natural isomorphisms $E\ten V \iso V
\liso V\ten E$ for each augmented $k$-quiver $V$. For two
augmented $k$-quivers $V$, $W$, we define $\Homq(V,W)$ to be the
$k$-quiver whose objects are the morphisms of augmented
$k$-quivers $f: V \ra W$ and such that for two morphisms of
augmented $k$-quivers $f$ and $g$, the $k$-module
$\Homq(V,W)(f,g)$ is formed by the families $(\phi_{x,y})$ of
\[
\prod_{x,y\in\obj(V)} \Hom_k(V(x,y), W(fx,gy))
\]
such that, for $f\neq g$, we have $\eps_{fx,gy} \phi_{x,y}=0$
for all $x,y$ of $V$, and for $f=g$, the map $\phi_{x,y}$ takes
the kernel of $\eps_{x,y}$ to the kernel of $\eps_{fx,fy}$ and 
the scalar
\[
\eps_{fx,fx} \circ \phi_{x,x} \circ \eta_{x,x}(\id_x)
\]
is independent of $x\in\obj(V)$. We endow $\Homq(V,W)$ with the
morphism $\eps$ which takes a family $(\phi_{x,y})$ to this scalar
and with the morphism $\eta$ which takes $\id_f$ to the family
\[
\phi_{x,y}: v \mapsto f(v).
\]
Clearly we obtain an augmented $k$-quiver $\Homq(V,W)$. For example,
we have $\Homq(E,W)\iso W$.

A {\em graded $k$-quiver} (resp. a {\em differential graded $k$-quiver})
is a $k$-quiver $V$ such that $V(x,y)$ is a $\Z$-graded $k$-module
(resp. a cochain complex of $k$-modules) for all objects $x,y$ of of $V$.
One defines the tensor product and the functor $\Homq$ for graded
and differential graded augmented $k$-quivers analogously.

\begin{lemma}The category of augmented $k$-quivers is monoidal with
unit object $E$ and $\Homq$ is an internal Hom-functor,
i.e. we have a functorial
bijection
\[
\Hom(U\ten V, W) \iso \Hom(U,\Homq(V,W))
\]
for all augmented $k$-quivers $U,V,W$. The same assertions hold for
the categories of graded and that of differential graded
augmented $k$-quivers.
\end{lemma}

\begin{proof} We define a map
\[
\Phi: \Hom(U\ten V, W) \iso \Hom(U,\Homq(V,W)).
\]
Let $F:U\ten V \to W$ be a morphism of augmented $k$-quivers. We
define
\[
\Phi F : U \to \Homq(V,W)
\]
as follows: For an object $u$ of $U$, we define $(\Phi F)(u)$ to
be the morphism of $k$-quivers $V\to W$ which sends an object $v$
of $V$ to $F(u,v)$ and a morphism $\beta: v_1 \to v_2$ of $V$ to
$F(\eta_u\ten \beta)$, where we write $\eta_u$ for $\eta_{u,u}$.
Then $(\Phi F)(u)$ is compatible with the unit and the
augmentation since $F(\eta_u\ten \eta_v)= \eta_{F(u,v)}$ and
\[
\eps_{F(u,v)} F (\eta_u \ten \beta) = (\eps_u \ten
\eps_v)(\eta_u\ten\beta)= \eps_v(\beta) \ko
\]
where we write $\eps_u$ instead of $\eps_{u,u}$. Now let $u_1$,
$u_2$ be objects of $U$ and $\alpha: u_1 \to u_2$ a morphism. We
define a morphism $(\Phi F)(\alpha): (\Phi F)(u_1) \to (\Phi
F)(u_2)$ of $\Homq(V,W)$ by sending a morphism $\beta :v_1 \to
v_2$ of $V$ to $F(\alpha\ten \beta): F(u_1,v_1)\to F(u_2,v_2)$. We
have
\[
\eps_{F(u_1,v_1), F(u_2,v_2)} (F(\alpha\ten\beta))=
\eps_{u_1,u_2}(\alpha) \eps_{v_1,v_2}(\beta)
\]
since $F$ is compatible with the augmentation. This shows that
$(\Phi F)(\alpha)$ lies indeed in $\Homq(V,W)$.

We define a map $\Psi : \Hom(U,\Homq(V,W)) \to \Hom(U\ten V, W)$.
Let $G$ be a morphism of augmented quivers from $U$ to
$\Homq(V,W)$. We define $\Psi G$ as follows: For objects $u$ of
$U$ and $v$ of $V$, we put $(\Psi G)(u,v)= (Gu)(v)$. For a
morphism $\alpha\ten \beta : (u_1,v_1) \to (u_2,v_2)$ of $U\ten
V$, we define
\[
(\Psi G)(\alpha\ten\beta) = (G \alpha)(\beta).
\]
Let us check that $\Psi G$ is a morphism of augmented $k$-quivers.
We have
\[
(\Psi G)(\eta_u\ten\eta_v) = (G \eta_u)(\eta_v).
\]
Since $G: U \to \Homq(V,W)$ respects $\eta$, the morphism
$G\eta_u$ equals $\eta_{Gu}$ and this is given by $\beta \mapsto
(Gu)(\beta)$. Since $Gu$ respects $\eta$, we have
$(Gu)(\eta_v)=\eta_{(Gu)(v)}$. Finally, we do obtain
\[
(\Psi G)(\eta_u\ten\eta_v) = \eta_{(Gu)(v)}.
\]
Now, for morphisms $\alpha:u_1\to u_2$ and $\beta:v_1\to v_2$, 
we compute
\begin{align*}
\eps((\Psi G)(\alpha\ten\beta)) &=
\eps((G \alpha)(\beta)) \\
&= \eps((G \alpha)(\eta \eps(\beta))) = 
\eps((G\alpha)(\eta(1))) \eps(\beta) \ko
\end{align*}
where the second last equality holds since $G\alpha$ vanishes on
the kernel of $\eps$ and we have omitted the subscripts for better
readability. Note that the scalar $\eps((G\alpha)(\eta(1)))$ is
the image of $G\alpha$ under the augmentation of $\Homq(V,W)$.
Since $G$ respects the augmentation, this scalar equals $\eps(\alpha)$.
So we have
\[
\eps((\Psi G)(\alpha\ten\beta)) = \eps (G\alpha) \eps(\beta)
=\eps(\alpha)\eps(\beta).
\]
Thus $\Psi G$ respects the augmentation. We have shown that the
maps $\Phi$ and $\Psi$ are well-defined. Clearly they are inverse
to each other.
\end{proof}

For an augmented $k$-quiver $V$, we denote by $\ol{V}$ the {\em reduction}
of $V$, i.e. the $k$-quiver with the same objects, and with
\[
\ol{V}(x,y) = \left\{ \begin{array}{ll} V(x,y) & \mbox{if } x \neq y \\
                                        V(x,x)/\im \eta_{x,x} & \mbox{if } x=y.
                      \end{array} \right.
\]
The reduction functor is an equivalence from the category of augmented
$k$-quivers to that of $k$-quivers but this equivalence is not monoidal
for the naive tensor product of $k$-quivers.

The augmented $k$-quivers with one object form a full monoidal subcategory of
the category of augmented $k$-quivers but this subcategory does not admit
an internal Hom-functor.

\subsection{Augmented cocomplete cocategories}

A {\em cocategory} is a $k$-quiver $C$ endowed with $k$-linear maps
\[
\Delta : C(x,y) \ra \bigoplus_{z\in \obj(C)} C(z,y)\ten C(x,z)
\]
for all objects $x,y$ of $C$ such that the natural coassociativity
condition holds.
It is {\em cocomplete} if each $f\in C(x,y)$ lies in the kernel of a
sufficiently high iterate of $\Delta$. It is {\em counital} if it is endowed
with a morphism of $k$-quivers $\eta : C \ra k\obj{C}$ such that the
two counit equations hold. An {\em augmented cocategory} is
a counital cocategory $C$ endowed with a morphism of counital cocategories
$\eps : k\obj(C) \ra C$ such that the composition of $\eta$ with $\eps$ is
the identity of $k\obj(C)$. Then it is in particular an augmented $k$-quiver
and the reduced quiver $\ol{C}$ becomes a cocategory. The augmented cocategory
$C$ is {\em cocomplete} if its reduction $\ol{C}$ is cocomplete.
These definitions admit obvious {\em graded and differential graded variants}.

Clearly the unit object $E$ of the category of augmented
$k$-quivers admits a unique structure of (cocomplete) augmented
cocategory. The tensor product $C_1\ten C_2$ of the $k$-quivers
underlying two cocomplete augmented cocategories admits
a natural structure of cocomplete augmented cocategory.
For this tensor product, the category of cocomplete augmented
cocategories becomes a monoidal category with unit object $E$.

Let $V$ be an augmented $k$-quiver. The {\em tensor cocategory $TV$} has
the same objects as $V$; and for any two objects $x,y$ of $TV$, the
$k$-module $(TV)\,(x,y)$ is the direct sum of $V(x,y)$ with the tensor products
\[
\ol{V}(z_n,y)\ten \ol{V}(z_{n-1},z_n) \ten \cdots \ten \ol{V}(x,z_1)
\]
where $n\geq 1$ and $(z_1, \ldots, z_n)$ runs through all sequences
of objects of $V$. For an element $(f_{n}, \ldots ,f_1)$, $n\geq 1$, of
$\ol{TV}(x,y)$, one puts
\begin{align*}
\Delta (f_{n},\ldots, f_1) = &   \eta_{y,y}\ten (f_{n},\ldots, f_1) +
\sum_{i=1}^{n} (f_{n}, \ldots, f_{i+1})\ten (f_i, \ldots, f_1)  \\
                               & +  (f_{n},\ldots, f_1)\ten \eta_{x,x}.
\end{align*}
The counit (resp.~the coaugmentation) of $TV$ are obtained by
composing the corresponding morphisms $V \ra k\obj(V)$ (resp. $k\obj(V)\ra V$)
with the projection $pr_1: TV\ra V$ (resp.~the inclusion $V\ra TV$).
The construction carries over to {\em augmented graded quivers} and
{\em augmented differential graded quivers}.

Denote by $\aqu$ the category of augmented $k$-quivers and by
$\aco$ that of cocomplete augmented $k$-cocategories.
The forgetful functor $F : \aco \ra \aqu$ is monoidal.

\begin{lemma}The forgetful functor admits the right adjoint $V \mapsto TV$.
More precisely, for each augmented cocategory $C$ and each augmented
$k$-quiver $V$, the map
\[
\Hom_{\aco}(C,TV) \ra \Hom_\aqu (FC,V) \ko f \mapsto pr_1 \circ f
\]
is bijective. Its inverse takes a morphism $g$ to the unique morphism
of augmented cocategories $G: C\ra TV$ whose reduction sends $c\in \ol{C}$
to
\[
\sum_{n=1}^\infty g^{\ten n} \Delta^{(n)}(c).
\]
The analogous assertions hold in the graded and in the differential
graded setup.
\end{lemma}

\begin{theorem}\label{closed-monoidal-thm}
The monoidal category of cocomplete augmented cocategories
admits an internal Hom-functor $\Homc$. Moreover, if $C$ is a cocomplete
augmented cocategory, then $\Homc(C,TV) = T\Homq(FC,V)$ for each
augmented dg $k$-quiver $V$.  The analogous assertions hold in the
graded and in the differential graded setting.
\end{theorem}

\begin{proof}We have to show that the functor $?\ten C$ admits a right
adjoint, i.e. that the functor $\Hom(?\ten C, C')$ is representable
for each $C'$. If $C'=TV$, then
\begin{eqnarray*}
\Hom(?\ten C,C') & = & \Hom(?\ten C, TV) = \Hom( F? \ten FC, V)  \\
                 & = & \Hom(F?, \Homq(FC,V)) = \Hom(?,T\Homq(FC,V)).
\end{eqnarray*}
If $C'$ is arbitrary, we can write it as the equalizer of a pair of
morphisms between cofree objects:
\[
TV_1 \stackrel{\ra}{\ra} TV_2.
\]
Then we obtain the representing object $\Homc(C,C')$ as the equalizer
of the pair
\[
\Homc(?,TV_1) \stackrel{\ra}{\ra} \Homc(?, TV_2).
\]
\end{proof}

\subsection{Explicit description of $\Homc(C,TV)$ and the action morphism}
Let $C$ be a cocomplete augmented cocategory and $V$ an augmented $k$-quiver.
For $D=TV$, we will describe $\Homc(C,D)$ and the adjunction
morphism (=action morphism)
\[
\Phi: \Homc(C,D)\ten C \to D.
\]
These generalize V.~Lyubashenko's constructions
\cite{Lyubashenko03}.
The objects of
\[
\Homc(C,TV) = T \Homc(FC,V).
\]
are the same as those of $\Homq(FC,V)$. Thus they are
in bijection with the morphisms of augmented $k$-quivers $f$
from $C$ to $V$. The morphism $\eta$ of $T\Homq(FC,V)$ is obtained
by composing that of $\Homq(FC,V)$ with the canonical inclusion and dually,
the morphism $\eps$ by composing that of $\Homq(FC,V)$ with the
projection. If $f$ and $g$ are two objects, a morphism from $f$
to $g$ in $T\Homq(FC,V)$ is given by a tensor product of $n\geq 0$
composable morphisms of $\Homq(FC,V)$. For $n=0$, this means that
$f=g$ and that the morphism is a multiple of $\eta_{f,f}(1)=f$;
for $n\geq 1$, the morphism is a tensor
\[
(\phi_n, \ldots, \phi_1)= \phi_n\ten \cdots \ten \phi_1
\]
of morphisms $\phi_i$ from $f_{i-1}$ to $f_i$ in the reduction
$\overline{\Homq(FC,V)}$, where $f_0=f$ and $f_n=g$.
Thus, $\phi_i$ belongs to
\[
\prod_{x,y\in \obj(C)} \Hom_k(C(x,y), V(f_{i-1}(x), f_i(y))
\]
and the composition
\[
\eps_{f_{i-1}(x), f_i(y)} \circ (\phi_i)_{x,y}
\]
vanishes for all $x,y$.
The image under $\Delta$ of such a
morphism $(\phi_n, \ldots, \phi_1)$ is
\[
\sum_{i=0}^{n} (\phi_n, \ldots, \phi_{i+1}) \ten (\phi_{i},
\ldots, \phi_1).
\]
Here, the terms for $i=0$ and $i=n$ are respectively
\[
(\phi_n, \ldots, \phi_1)\ten f \mbox{ and } g \ten (\phi_n,
\ldots, \phi_1) \ko
\]
by definition of the morphism $\eta$ for $\Homq(FC,V)$. When we
apply an iterate $\Delta^{(N)}$ of $\Delta$ to $(\phi_n, \ldots,
\phi_1)$, we  obtain sums of tensor products of $N$ factors
each of which may be reduced to some $f_i$. For example, we get,
among others, the
term
\[
f_7 \ten (\phi_7) \ten f_6 \ten (\phi_6,\phi_5) \ten f_4 \ten f_4
\ten (\phi_4, \phi_3, \phi_2) \ten f_1 \ten (\phi_1)
\]
for $n=7$ and $N=9$.

The adjunction morphism $\Phi: \Homc(C,TV)\ten C \to TV$ is the
unique morphism of cocomplete augmented cocategories whose image
under the forgetful functor makes the square
\[
\xymatrix{ T\Homq(FC,V) \ten C \ar[r] \ar[d] & TV \ar[d] \\
\Homq(FC,V)\ten C \ar[r] & V }
\]
commutative in the category of augmented $k$-quivers. Thus, the
morphism $\Phi$ takes an object $(f,x)$ of the tensor product
$T\Homq(FC,V)\ten C$ to $f(x)$ and a morphism $(\phi_n, \ldots,
\phi_1)\ten c$
to the image in
$\bigoplus_{N\geq n} \ol{V}\tp{N}$ of the sum of the
\[
(f_n\tp{r_n} \ten \phi_n \ten f_{n-1}\tp{r_{n-1}} \ten \cdots
\ten f_1\tp{r_1} \ten \phi_1 \ten f_0\tp{r_0}) \circ \Delta^{(N)}(c)
\]
where $N=n+\sum r_i$. Note that this sum has indeed only finitely many
non zero terms: Since $C$ is cocomplete, $\Delta^{(M)}(c)$
vanishes in $\overline{C}\tp{M}$ for some $M\gg 0$. Therefore,
$\Delta^{(N)}(c)$ vanishes in $C\tp{p}\ten \overline{C}\tp{M} \ten
C\tp{q}$ for all $N=p+M+q$, $p,q\geq 0$. Now for large enough $N$,
each term of the image of $(\phi_n, \ldots, \phi_1)$ under
$\Delta^{(N)}$ will contain a chain of $M$ consecutive $f_i$'s and therefore
induce a map from $C\tp{p}\ten \overline{C}\tp{M} \ten
C\tp{q}$ to $V\tp{p}\ten \overline{V} \tp{M}\ten V\tp{q}$. Hence
the image of
\[
\Delta^{(N)}(\phi_n, \ldots, \phi_1)\circ \Delta^{(N)}(c)
\]
in $\ol{V}\tp{N}$ vanishes.

Note that if $C$ is the unit object $E$, the action morphism
$\Homc(E,D)\ten E \to D$ is the identity. More generally, if
$C$ is non empty, then $E$ is a retract of $C$ and the action
morphism $\Homc(C,D)\ten C\to D$ admits a section. Dually,
the coaction morphism $D\to\Homc(C,D\ten C)$ admits a retraction
in this case.

\subsection{The composition morphism} \label{composition-morphism}
Let $C$ be a cocomplete
augmented cocategory and $V$, $W$ augmented $k$-quivers. Since
the functor $?\ten C$ is left adjoint to $\Homc(?,C)$, there is a unique
morphism, {\em the composition morphism},
\[
\zeta: \Homc(TV,TW)\ten \Homc(C,TV) \longrightarrow \Homc(C,TW)
\]
which yields a commutative square
\[
\xymatrix{
\Homc(TV,TW) \ten \Homc(C,TV) \ten C \ar[d]_{\id\ten\Phi}
\ar[r]^-{\zeta \ten C} &
\Homc(C,TW)\ten C \ar[d]^{\Phi} \\
\Homc(TV,TW) \ten TV  \ar[r]^-{\Phi} & TW.}
\]
The composition morphism is constructed as the unique morphism
such that the following square commutes in the category of
augmented $k$-quivers
\[
\xymatrix{
T\Homq(TV,W)\ten T\Homq(C,V) \ar[r] \ar[d] & T\Homq(C,W) \ar[d]\\
\Homq(TV,W)\ten T\Homq(C,V) \ar[r] & \Homq(C,W),
}
\]
where the bottom arrow is described as follows: a pair
$(f,f')$ of morphisms is sent to the composition $f\circ f'$.
A morphism $\phi'\ten (\phi_n, \ldots, \phi_1)$ is sent to
the morphism
\[
c \mapsto \sum \phi'\circ (f_n\tp{r_n}\ten \phi_n \ten \cdots
\ten f_i\tp{r_i}\ten \phi_i \ten f_{i-1}\tp{r_{i-1}} \ten \cdots
\ten \phi_1 \ten f_0\tp{r_0}) \circ \Delta^{(N)} (c)
\]
where $N$ is the number of arguments of $\phi'$ and
$n+\sum r_i=N$.

\subsection{The comparison morphism} Let $C$ be a cocomplete
augmented cocategory and $W$ an augmented $k$-quiver. The
{\em comparison morphism} is the composition
\[
F\Homc(C,TW) \to \Homq(FC,F\Homc(C,TW)\ten FC) \to \Homq(FC,FTW)
\]
of the coaction morphism for $?\ten FC$  with
the action morphism for $?\ten C$.
This morphism induces the natural injection in the sets
of objects.  It takes a morphism $(\phi_n, \ldots, \phi_1)$, $n\geq
1$, of $T\Homq(C,W)$ to the morphism of $\Homq(FC,FTW)$ which sends
$c\in C(x,y)$ to the image in $\bigoplus_{N\geq n} \ol{W}\tp{N}$ of
the sum of the
\[
(f_n\tp{r_n} \ten \phi_n \ten f_{n-1}\tp{r_{n-1}} \ten \cdots
\ten f_1\tp{r_1} \ten \phi_1 \ten f_0\tp{r_0}) \circ \Delta^{(N)}(c)
\]
where $N=n+\sum r_i$. We have seen above that this sum has only
finitely many non zero terms. If the map
$C \to \overline{C}\tp{n}$ induced by the iterated comultiplication
is surjective for all $n\geq 2$, then
the comparison morphism is clearly injective on the sets
of morphisms. This happens for example if $C=TV$ for some
augmented $k$-quiver $V$. The comparison morphism is
compatible with compositions.

\subsection{The differential graded case}
Let $\adgco$ denote the category of cocomplete augmented dg coalgebras and
$\agco$ the category of cocomplete augmented graded coalgebras.
The functor
\[
\For :\adgco \ra \agco
\]
which forgets the differentials is monoidal. Thus, if $A,B$ are
objects of $\adgco$, the adjunction morphism
\[
\Homq(A,B)\ten A \to B
\]
yields a morphism
\[
\For\Homq(A,B)\ten \For A \to \For B
\]
and hence a morphism
\[
\For \Homq(A,B) \ra \Homq(\For A,\For B).
\]

\begin{lemma} The above morphism induces isomorphisms in the morphism
spaces, \ie $\For \Homq(A,B)$ is the full subcocategory of
$\Homq(\For A,\For B)$ whose objects are the morphisms of $\agco$
which commute with the differential.
\end{lemma}

\begin{proof} It suffices to prove this when $B$ is cofree on some
augmented differential graded $k$-quiver $V$. Then we have
\[
\For \Homq(A, TV) = \For T \Homq(A, V) = T \For \Homq(A,V)
\]
where $\Homq(A,V)$ is the internal $\Homq$-object of the category of
augmented differential graded $k$-quivers and we have also denoted by
$\For$ the forgetful functor from dg $k$-quivers to graded $k$-quivers.
Now it is clear that we have a full morphism of graded quivers
\[
\For\Homq(A,V) \ra \Homq(\For A,\For V)
\]
and this implies the assertion.
\end{proof}

\subsection{$A_\infty$-functor categories} \label{Ainf-functors}
For an augmented graded $k$-quiver $V$, denote by $SV$ the augmented
graded $k$-quiver with the same objects whose reduction has the
shifted $k$-modules $S \overline{V}(x,y)$ as morphism spaces.  An {\em
augmented $A_\infty$-category} is an augmented graded $k$-quiver $A$
whose bar construction $BA=TSA$ is endowed with a graded coalgebra
differential so that it becomes a cocomplete augmented dg
cocategory. Equivalently, the reduction $\overline{A}$ is endowed with
a structure of $A_\infty$-category without units in the sense of
Kontsevich-Soibelman \cite[4.1]{KontsevichSoibelman01}. For example,
if $A=k\oplus \overline{A}$ is an ordinary augmented associative
algebra, then we view it as an $A_\infty$-category with one
object by endowing
\[
BA=TSA= k\oplus \bigoplus_{n\geq 1} (S\ol{A})\tp{n}
\]
with the unique coalgebra differential such that $d(x\ten y) = xy$
for $x,y\in\ol{A}$. Then the underlying complex of $BA$ is
\[
k \la \ol{A} \la \ol{A}\ten \ol{A} \la \cdots
\]
and $H_*(BA)= \Tor_*^A(k,k)$.

Let $A$, $A'$ be two augmented $A_\infty$-categories.
By theorem~\ref{closed-monoidal-thm}, the cocategory
$\Homc(BA, BA')$ is (canonically) cofree.
We define the {\em $A_\infty$-functor category $\Funinf(A,A')$}
by
\[
\Homc(BA,BA') = B(\Funinf(A,A')).
\]
It is an augmented $A_\infty$-category.
By definition, its objects are the {\em $A_\infty$-functors} from $A$
to $A'$. According to theorem~\ref{closed-monoidal-thm} and the above
lemma, the underlying graded category of $\Homc(BA,BA')$ is a full
subcategory of
\[
\Homc(BA, TSA')= T\Homq(BA,SA') \ko
\]
Thus, the underlying
graded quiver of $S\Funinf(A,A')$ is a full subquiver of
\[
\Homq(BA,SA').
\]
Suppose that $A$ and $A'$ are ordinary augmented algebras. Then the objects
of $\Funinf(A,A')$ identify with the morphisms of augmented dg coalgebras
$BA \to BA'$. For any two objects $f,g$, the morphism space
\[
S\Funinf(A,A')(f,g)
\]
is
\[
\Hom(k,S\ol{A'}) \oplus \Hom(\ol{BA}, S\ol{A'})
\]
for $f\neq g$ and
\[
 k \oplus \Hom(k,S\ol{A'})
                                       \oplus \Hom(\ol{BA}, S\ol{A'})
\]
for $f=g$. {\em If $f\neq g$}, then the morphisms $f\to g$ are
in canonical bijection with the space of {\em $(f,g)$-coderivations}
$D: BA \to BA'$, \ie the homogeneous linear maps satisfying
\[
\Delta \circ D = (f\ten D + D \ten g) \circ \Delta.
\]
Note that $\Funinf(A,A')$ has infinitely many objects, in general, even though
$A$ and $A'$ only have one. If $A=A'$ and $f=g$ is the identity
functor then
\[
\Funinf(A,A)(\id,\id)/k
\]
is the reduced Hochschild complex
\[
\ol{A} \to \Hom(\ol{A}, \ol{A}) \to \Hom(\ol{A} \ten \ol{A} , \ol{A})
\to \ldots
\to \Hom(\ol{A}\tp{p}, \ol{A}) \to \ldots
\]
where $\ol{A}$ appears in degree $-1$.

If $A,A', A''$ are three augmented $A_\infty$-categories,
the composition morphism
\[
\Homc(BA', BA'') \ten \Homc(BA,BA') \to \Homc(BA,BA'')
\]
yields a morphism
\[
B\Funinf(A',A'') \ten B\Funinf(A,A') \to B\Funinf(A,A'') \ko
\]
which we view as an `$A_\infty$-bifunctor'
\[
\Funinf(A',A'') \times \Funinf(A,A') \to \Funinf(A,A'').
\]
Let us show how this bifunctor yields the bialgebra structure on the
bar construction of the Hochschild complex.

We need the following
remark on cocategories: Let $\cc$ be an augmented $k$-cocategory
and $x$ and object of $\cc$. Then the endomorphism space
$\cc(x,x)$ is not, in general, a coalgebra in a natural way,
because the comultiplication
\[
\Delta: \cc(x,x) \to \bigoplus_y \cc(y,x) \ten \cc(x,y)
\]
does not take values in $\cc(x,x)\ten \cc(x,x)$. However, if we
define $\End(x)$ as the set of all $c\in \cc(x,x)$ such that
$\Delta^{(N)}(c)$ lies in $\cc(x,x)\tp{N}$ for all $N\geq 2$,
then $\End(x)$ becomes an augmented coalgebra, which we call
the {\em endomorphism coalgebra of $x$}.

Now let $A$ be an ordinary
augmented algebra and consider the cocategory
\[
\cc=\Homc(BA,BA) = B\Funinf(A,A).
\]
Then the endomorphism coalgebra of the identity functor
$\End(\id)$ is the bar construction $B(C(A,A))$, on the
reduced Hochschild complex $C(A,A)$, \ie
\[
\ol{A} \to \Hom(\ol{A}, \ol{A}) \to \Hom(\ol{A} \ten \ol{A} , \ol{A})
\to \ldots
 \to \Hom(\ol{A}\tp{p}, \ol{A}) \to \ldots
\]
with $\ol{A}$ appearing in degree $0$. The composition bifunctor
\[
B\Funinf(A,A)\ten B\Funinf(A,A)\to B\Funinf(A,A)
\]
then yields an associative multiplication on the
endomorphism coalgebra
\[
\End(\id)=B(C(A,A))
\]
so that it becomes a bialgebra. Using the explicit description
of the composition morphism in section~\ref{composition-morphism}
one checks that this bialgebra structure is given by the
brace operations \cite{Baues81} \cite{Kadeishvili88} \cite{Getzler93}
 \cite{GetzlerJones94}.

\subsection{Modules and bimodules} Let $\cc(k)$ denote the
dg category of complexes of vector spaces and $\cc(k)^+$
the augmented dg category whose reduction is $\cc(k)$.
Let $A$ be an augmented $A_\infty$-category. By definition,
the category of (left) $A_\infty$-modules over $A$ is the category
\[
\Funinf(A, \cc(k)^+).
\]
It is in fact a dg category (since the target category is
a dg category). If $A'$ is another augmented $A_\infty$-category,
then we can consider
\[
\Funinf(A', \Funinf(A,\cc(k)^+)).
\]
From the canonical isomorphism
\[
\Homc(U, \Homc(V,W)) \iso \Homc(U\ten V, W)
\]
in the category of cocomplete augmented dg cocategories,
we obtain
\[
\Funinf(A', \Funinf(A, \cc(k)^+) \iso \Homc(BA' \ten BA, \cc(k)^+),
\]
and the last term is the category of $A_\infty$-bimodules. Equivalences
of this type are quite unpleasant to prove directly,
\cf for example Lemme 5.3.0.1 of \cite{Lefevre03}.


\begin{thebibliography}{10}

\bibitem{Baues81}
H.~J. Baues, \emph{The double bar and cobar constructions}, Compositio Math.
  \textbf{43} (1981), no.~3, 331--341.

\bibitem{BeilinsonGinzburgSoergel96}
Alexander Beilinson, Victor Ginzburg, and Wolfgang Soergel, \emph{Koszul
  duality patterns in representation theory}, J. Amer. Math. Soc. \textbf{9}
  (1996), no.~2, 473--527.

\bibitem{BensonKrauseSchwede04}
David Benson, Henning Krause, and Stefan Schwede, \emph{Realizability of
  modules over {T}ate cohomology}, Trans. Amer. Math. Soc. \textbf{356} (2004),
  no.~9, 3621--3668 (electronic).

\bibitem{BespalovLyubashenkoManzyuk06}
Yuri Bespalov, V.~V. Lyubashenko, and Oleksandr Manzyuk, \emph{Closed
  multicategory of (pretriangulated) ${A}_\infty$-categories}, in progress,
  2006.

\bibitem{Beilinson78}
A.~A. Be\u{\i}linson, \emph{Coherent sheaves on {${\bf P}\sp{n}$} and problems
  in linear algebra}, Funktsional. Anal. i Prilozhen. \textbf{12} (1978),
  no.~3, 68--69.

\bibitem{BondalVandenBergh03}
A.~Bondal and M.~van~den Bergh, \emph{Generators and representability of
  functors in commutative and noncommutative geometry}, Mosc. Math. J.
  \textbf{3} (2003), no.~1, 1--36, 258.

\bibitem{ChuangKingxy}
Joseph Chuang and Alastair King, \emph{The ${A}_\infty$-structure on the
  {Y}oneda algebra of a monomial algebra}, preliminary notes.

\bibitem{DwyerSpalinski95}
W.~G. Dwyer and J.~Spali{\'n}ski, \emph{Homotopy theories and model
  categories}, Handbook of algebraic topology, North-Holland, Amsterdam, 1995,
  pp.~73--126.

\bibitem{Floystad00}
Gunnar Floystad, \emph{{Koszul duality and equivalences of categories}},
  arXiv:math.RA/0012264.

\bibitem{Fukaya93}
Kenji Fukaya, \emph{Morse homotopy, {$A\sp \infty$}-category, and {F}loer
  homologies}, Proceedings of GARC Workshop on Geometry and Topology '93
  (Seoul, 1993) (Seoul), Lecture Notes Ser., vol.~18, Seoul Nat. Univ., 1993,
  pp.~1--102.

\bibitem{Getzler93}
Ezra Getzler, \emph{Cartan homotopy formulas and the {G}auss-{M}anin connection
  in cyclic homology}, Quantum deformations of algebras and their
  representations (Ramat-Gan, 1991/1992; Rehovot, 1991/1992), Israel Math.
  Conf. Proc., vol.~7, Bar-Ilan Univ., Ramat Gan, 1993, pp.~65--78.

\bibitem{GetzlerJones94}
Ezra Getzler and J.~D.~S. Jones, \emph{Operads, homotopy algebra, and iterated
  integrals for double loop spaces}, hep-th/9403055.

\bibitem{GetzlerJones90}
Ezra Getzler and John D.~S. Jones, \emph{{$A\sb \infty$}-algebras and the
  cyclic bar complex}, Illinois J. Math. \textbf{34} (1990), no.~2, 256--283.

\bibitem{Hinich97}
Vladimir Hinich, \emph{Homological algebra of homotopy algebras}, Comm. Algebra
  \textbf{25} (1997), no.~10, 3291--3323.

\bibitem{Hovey99}
Mark Hovey, \emph{Model categories}, Mathematical Surveys and Monographs,
  vol.~63, American Mathematical Society, Providence, RI, 1999.

\bibitem{Huebschmann86}
Johannes Huebschmann, \emph{The homotopy type of {$F\Psi\sp q$}. {T}he complex
  and symplectic cases}, Applications of algebraic $K$-theory to algebraic
  geometry and number theory, Part I, II (Boulder, Colo., 1983), Contemp.
  Math., vol.~55, Amer. Math. Soc., Providence, RI, 1986, pp.~487--518.

\bibitem{Huebschmann89}
\bysame, \emph{The mod-{$p$} cohomology rings of metacyclic groups}, J. Pure
  Appl. Algebra \textbf{60} (1989), no.~1, 53--103.

\bibitem{Kadeishvili88}
T.~V. Kadeishvili, \emph{The structure of the {$A(\infty)$}-algebra, and the
  {H}ochschild and {H}arrison cohomologies}, Trudy Tbiliss. Mat. Inst. Razmadze
  Akad. Nauk Gruzin. SSR \textbf{91} (1988), 19--27.

\bibitem{Kadeishvili80}
T.~V. Kadei{\v{s}}vili, \emph{On the theory of homology of fiber spaces},
  Uspekhi Mat. Nauk \textbf{35} (1980), no.~3(213), 183--188, International
  Topology Conference (Moscow State Univ., Moscow, 1979).

\bibitem{Keller01}
Bernhard Keller, \emph{Introduction to {$A$}-infinity algebras and modules},
  Homology Homotopy Appl. \textbf{3} (2001), no.~1, 1--35 (electronic).

\bibitem{Keller02}
\bysame, \emph{Addendum to: ``{I}ntroduction to {$A$}-infinity algebras and
  modules'' [{H}omology {H}omotopy {A}ppl.\ {\bf 3} (2001), no.\ 1, 1--35;},
  Homology Homotopy Appl. \textbf{4} (2002), no.~1, 25--28 (electronic).

\bibitem{Kontsevich94}
Maxim Kontsevich, \emph{Homological algebra of mirror symmetry}, Proceedings of
  the International Congress of Mathematicians, Vol.\ 1, 2 (Z\"urich, 1994)
  (Basel), Birkh\"auser, 1995, pp.~120--139.

\bibitem{Kontsevich98}
\bysame, \emph{Triangulated categories and geometry}, Course at the {\'{E}}cole
  {N}ormale Sup\'erieure, Paris, Notes taken by J. Bella\"{\i}che, J.-F. Dat,
  I. Marin, G. Racinet and H. Randriambololona, 1998.

\bibitem{KontsevichSoibelman00}
Maxim Kontsevich and Yan Soibelman, \emph{Deformations of algebras over operads
  and the {D}eligne conjecture}, Conf{\'e}rence Mosh{\'e} Flato 1999, Vol. I
  (Dijon), Math. Phys. Stud., vol.~21, Kluwer Acad. Publ., Dordrecht, 2000,
  pp.~255--307.

\bibitem{KontsevichSoibelman01}
\bysame, \emph{Homological mirror symmetry and torus fibrations}, Symplectic
  geometry and mirror symmetry (Seoul, 2000), World Sci. Publishing, River
  Edge, NJ, 2001, pp.~203--263.

\bibitem{Lefevre03}
Kenji Lef\`evre-Hasegawa, \emph{Sur les ${A}_\infty$-cat\'egories}, Th\`ese de
  doctorat, {U}niversit\'e {D}enis {D}iderot -- Paris 7, November 2003,
  available at B.~Keller's homepage.

\bibitem{LuPalmieriWuZhang04}
D.~M. Lu, J.~H. Palmieri, Q.~S. Wu, and J.~J. Zhang, \emph{{$A\sb
  \infty$}-algebras for ring theorists}, Proceedings of the International
  Conference on Algebra, vol.~11, 2004, pp.~91--128.

\bibitem{Lyubashenko03}
Volodymyr Lyubashenko, \emph{Category of {$A\sb \infty$}-categories}, Homology
  Homotopy Appl. \textbf{5} (2003), no.~1, 1--48 (electronic).

\bibitem{Madsen02}
Dag Madsen, \emph{Ph. d. thesis}, NTNU, Trondheim, 2002.

\bibitem{May69}
J.~Peter May, \emph{Matric {M}assey products}, J. Algebra \textbf{12} (1969),
  533--568.

\bibitem{Mitchell72}
Barry Mitchell, \emph{Rings with several objects}, Advances in Math. \textbf{8}
  (1972), 1--161.

\bibitem{PenkavaSchwarz95}
Michael Penkava and Albert Schwarz, \emph{{$A\sb \infty$} algebras and the
  cohomology of moduli spaces}, Lie groups and Lie algebras: E. B. Dynkin's
  Seminar, Amer. Math. Soc. Transl. Ser. 2, vol. 169, Amer. Math. Soc.,
  Providence, RI, 1995, pp.~91--107.

\bibitem{Proute84}
Alain Prout\'e, \emph{Alg\`ebres diff\'erentielles fortement homotopiquement
  associatives}, Th\`ese d'\'Etat, Universit\'e Paris VII, 1984.

\bibitem{Quillen67}
Daniel~G. Quillen, \emph{Homotopical algebra}, Lecture Notes in Mathematics,
  No. 43, Springer-Verlag, Berlin, 1967.

\bibitem{Rickard89}
Jeremy Rickard, \emph{Morita theory for derived categories}, J. London Math.
  Soc. \textbf{39} (1989), 436--456.

\bibitem{Rouquier03}
Rapha\"el Rouquier, \emph{Dimensions of triangulated categories},
  arXiv:math.CT/0310134.

\bibitem{Seidel02}
Paul Seidel, \emph{Fukaya categories and deformations}, Proceedings of the
  International Congress of Mathematicians, Vol. II (Beijing, 2002) (Beijing),
  Higher Ed. Press, 2002, pp.~351--360.

\bibitem{Smirnov80}
V.~A. Smirnov, \emph{Homology of fiber spaces}, Uspekhi Mat. Nauk \textbf{35}
  (1980), no.~3(213), 227--230, International Topology Conference (Moscow State
  Univ., Moscow, 1979).

\bibitem{Stasheff63}
James~Dillon Stasheff, \emph{Homotopy associativity of {$H$}-spaces. {I},
  {II}}, Trans. Amer. Math. Soc. 108 (1963), 275-292; ibid. \textbf{108}
  (1963), 293--312.

\bibitem{Stasheff92}
Jim Stasheff, \emph{Differential graded {L}ie algebras, quasi-{H}opf algebras
  and higher homotopy algebras}, Quantum groups (Leningrad, 1990), Lecture
  Notes in Math., vol. 1510, Springer, Berlin, 1992, pp.~120--137.

\end{thebibliography}

\def\cprime{$'$}
\providecommand{\bysame}{\leavevmode\hbox to3em{\hrulefill}\thinspace}
\providecommand{\MR}{\relax\ifhmode\unskip\space\fi MR }
\providecommand{\MRhref}[2]{%
  \href{http://www.ams.org/mathscinet-getitem?mr=#1}{#2}
}
\providecommand{\href}[2]{#2}

\end{document}